\documentclass[10pt,fleqn]{amsproc}

\usepackage[utf8]{inputenc}
\usepackage[T1]{fontenc}
\usepackage{amsmath,amssymb,amsthm}
\usepackage{lmodern}

\usepackage{mathtools}
\usepackage{mathrsfs}
\usepackage{graphicx}

\usepackage{tikz}
\usepackage{tikz-cd}

\usetikzlibrary{arrows.meta,calc,decorations.markings,intersections}
\usepackage{xcolor}
\usepackage{hyperref}
\hypersetup{colorlinks=true,linkcolor=magenta,citecolor=blue,urlcolor=blue}

\usepackage{enumitem}

\newcommand{\textupbf}[1]{\textbf{\textup{#1}}}

\allowdisplaybreaks

\SetEnumitemKey{auto}{label=\textupbf{(\arabic*)},ref=(\arabic*),leftmargin=*}
\SetEnumitemKey{autoarabic}{label=\textupbf{(\arabic*)},ref=(\arabic*),leftmargin=*}
\SetEnumitemKey{autoalpha}{label=\textupbf{(\alph*)},ref=(\alph*),leftmargin=*}
\SetEnumitemKey{autoroman}{label=\textupbf{(\roman*)},ref=(\roman*),leftmargin=*}
\SetEnumitemKey{autoAlpha}{label=\textupbf{(\Alph*)},ref=(\Alph*),leftmargin=*}
\SetEnumitemKey{autoRoman}{label=\textupbf{(\Roman*)},ref=(\Roman*),leftmargin=*}

\newcommand{\rmop}[1]{\operatorname{#1}}

\newcommand{\frakop}[1]{\operatorname{\mathfrak{#1}}}
\newcommand{\calop}[1]{\operatorname{\mathcal{#1}}}

\renewcommand{\emph}[1]{\textit{\textbf{\textcolor{blue!80!black}{#1}}}}

\numberwithin{equation}{section}
\numberwithin{figure}{section}

\newtagform{bold}{\bfseries(}{\bfseries)}
\usetagform{bold}

\newcommand{\purplebfeq}[1]{\textcolor{purple!80!black}{\bfseries(#1)}}
\newtagform{bfpurple}[\purplebfeq]{}{}
\usetagform{bfpurple}

\newcommand{\frakU}{\mathfrak{U}}

\newcommand{\frakc}{\mathfrak{c}}

\newcommand{\frakg}{\mathfrak{g}}

\newcommand{\bbA}{\mathbb{A}}

\newcommand{\bbC}{\mathbb{C}}

\newcommand{\bbF}{\mathbb{F}}
\newcommand{\bbG}{\mathbb{G}}

\newcommand{\bbP}{\mathbb{P}}
\newcommand{\bbQ}{\mathbb{Q}}
\newcommand{\bbR}{\mathbb{R}}

\newcommand{\bbZ}{\mathbb{Z}}

\newcommand{\calB}{\mathcal{B}}
\newcommand{\calC}{\mathcal{C}}
\newcommand{\calD}{\mathcal{D}}

\newcommand{\calN}{\mathcal{N}}
\newcommand{\calO}{\mathcal{O}}

\newcommand{\calR}{\mathcal{R}}
\newcommand{\calS}{\mathcal{S}}

\newcommand{\calW}{\mathcal{W}}

\newcommand{\scrC}{\mathscr{C}}
\newcommand{\scrD}{\mathscr{D}}

\newcommand{\scrT}{\mathscr{T}}

\newcommand{\scrX}{\mathscr{X}}
\newcommand{\scrY}{\mathscr{Y}}

\newcommand{\bfG}{\mathbf{G}}
\newcommand{\bfH}{\mathbf{H}}

\newcommand{\bfJ}{\mathbf{J}}

\newcommand{\bfO}{\mathbf{O}}

\newcommand{\bfS}{\mathbf{S}}
\newcommand{\bfT}{\mathbf{T}}

\newcommand{\bfm}{\mathbf{m}}

\newcommand{\rmf}{\mathrm{f}}

\newcommand{\rmm}{\mathrm{m}}

\usepackage[backend=biber,style=alphabetic,maxbibnames=99,maxalphanames=5,
  giveninits=true,doi=true,isbn=false,url=true]{biblatex}
\AtBeginBibliography{\renewcommand{\emph}{\textit}}

\providecommand{\thmabovespace}{6pt}
\providecommand{\thmbelowspace}{6pt}

\newtheoremstyle{standardblue}%
  {\thmabovespace}{\thmbelowspace}
  {} 
  {0pt}{\bfseries\color{blue!80!black}}{.}{1em}
  {\thmname{#1}\thmnumber{~#2}\thmnote{~({\normalfont{#3}})}}
  
\newtheoremstyle{standarditred}%
  {\thmabovespace}{\thmbelowspace}
  {\itshape} 
  {0pt}{\bfseries\color{red!80!black}}{.}{1em}
  {\thmname{#1}\thmnumber{~#2}\thmnote{~({\normalfont{#3}})}}

\theoremstyle{standarditred}

\newtheorem{thm}{Theorem}

\newtheorem{proposition}{Proposition}[section]
\newtheorem{theorem}[proposition]{Theorem}
\newtheorem{lemma}[proposition]{Lemma}

\newtheorem*{claim*}{Claim}

\theoremstyle{standardblue}

\newtheorem*{remark*}{Remark}

\DeclareMathOperator{\nrd}{nrd}

\DeclareMathOperator{\trd}{trd}

\DeclareMathOperator{\End}{End}
\DeclareMathOperator{\cl}{cl}

\newcommand{\xto}{\xrightarrow}
\newcommand{\medsim}{\ensuremath\mathrel{\scalebox{1.05}{\ensuremath\sim}}}
\newcommand{\stosim}[1][]{\stackrel{\medsim}{\smash{\xrightarrow[#1]{\hspace{0.5em}}}\rule{0pt}{0.2ex}}}

\newcommand{\m}[2][p]{\begin{#1matrix}#2\end{#1matrix}}

\tikzset{
  symbol/.style={
    draw=none,
    every to/.append style={
      edge node={node [sloped, allow upside down, auto=false]{\(#1\)}}}
  }
}

\title[Equidistribution of Gross lattices]{Equidistribution of Gross lattices of supersingular elliptic curves}
\author{Lam L. Pham}
\email{lam.pham@ulb.be}
\address{Department of Computer Science, Universit\'e Libre de Bruxelles}
\thanks{L.L. Pham (ULB) is supported by the FRS-FNRS through grant number T.0238.25}

\newcommand{\Z}{\bbZ}
\newcommand{\R}{\bbR}
\newcommand{\Q}{\bbQ}

\newcommand{\C}{\bbC}
\DeclareMathOperator{\SL}{SL}
\DeclareMathOperator{\bfSL}{\mathbf{SL}}
\DeclareMathOperator{\GL}{GL}
\DeclareMathOperator{\bfGL}{\mathbf{GL}}

\DeclareMathOperator{\bfPGL}{\mathbf{PGL}}
\DeclareMathOperator{\SO}{SO}
\DeclareMathOperator{\bfSO}{\mathbf{SO}}
\DeclareMathOperator{\rmSS}{SS}
\DeclareMathOperator{\Aut}{Aut}
\DeclareMathOperator{\Ad}{Ad}
\DeclareMathOperator{\pr}{pr}
\DeclareMathOperator{\sh}{sh}
\DeclareMathOperator{\rmht}{ht}
\DeclareMathOperator{\Av}{Av}
\DeclareMathOperator{\frsl}{\frakop{sl}}

\newcommand{\mass}{\mathrm{mass}}
\newcommand{\bs}{\backslash}

\newtheorem{definition}[proposition]{Definition}
\newtheorem*{definition*}{Definition}

\newcommand{\onto}{\twoheadrightarrow}
\newcommand{\into}{\hookrightarrow}

\newcommand{\bfmu}{\boldsymbol{\mu}}
\newcommand{\bflambda}{\boldsymbol{\lambda}}
\renewcommand{\bar}{\overline}

\begin{document}

\begin{abstract}
For a prime \(p\), we consider two important families of integral quadratic lattices associated to isomorphism classes of supersingular elliptic curves over \(\bar{\bbF_{p}}\), the lattices of endomorphisms of trace zero, and the Gross lattices of endomorphisms.
We prove effective equidistribution of the shapes of these lattices as \(p\to\infty\) through primes, in the arithmetic locally symmetric space of unimodular lattices up to isometry.
As an application, we obtain, as \(p\to\infty\), the distribution of successive minima of these lattices.
These results are proved using techniques from homogeneous dynamics.
\end{abstract}

\maketitle

\tableofcontents

\section{Introduction and main result}

\subsection{Setting and notations}
\label{section:setting}

\subsection*{Elliptic curves}

Let \(p\geq 5\) be a prime and let \(\rmop{SS}_{p}\) denote the finite set of geometric isomorphism classes over \(\bar{\bbF}_{p}\) of supersingular elliptic curves, with \(h_{p}=|\rmop{SS}_{p}|\).
Let \(E\) be such a curve.
By Deuring's theorem, \(\calO_{E}=\End_{\bar{\bbF}_{p}}(E)\) is  a maximal order in the rational quaternion algebra \(B_{p,\infty}=\calO_{E}\otimes\Q\) ramified at \(p\) and \(\infty\).
For an order \(\calO\subset B_{p,\infty}\), we define the \emph{trace-zero lattice} \(\calO^{0}\) and the \emph{Gross lattice} \(T(\calO)\) as
\[
\calO^{0}=\ker(\trd:\calO\to \Z),\quad
T(\calO)=\{2x-\trd(x)\mid x\in \calO\}.
\]
In fact, we have \(T(\calO)=\calO^{0}\cap (\Z+2\calO)\).
Let \([E]\in \rmop{SS}_{p}\).
For the order \(\calO_{E}\), we denote these lattices by \(L_{E}\) and \(\Lambda_{E}\) respectively, with the quadratic form \(q_{E}\):
\begin{equation}
\label{eq:lattices}
L_{E} = \ker(\trd:\calO_{E} \to \bbZ),\quad
\Lambda_{E} = (\Z+2\calO_{E})\cap L_{E},\quad
q_{E} = \nrd|_{L_{E}}.
\end{equation}
Both lattices have rank \(3\), we have \(2L_{E}\subset \Lambda_{E}\subset L_{E}\), \([L_{E}:\Lambda_{E}]=4\), and \(q_{E}\) is a positive definite quadratic form on both \(L_{E}\) and \(\Lambda_{E}\).

There has been increasing interest in Gross lattices of supersingular elliptic curves in recent years, see for instance, \cite{GorenLove2025,HeKorpalTranVincent2026,AubryOyonoVincent2026}, and the references therein.

\subsection*{Shapes}

We now define our moduli space of unimodular lattices.
Let
\[
G_{\infty} = \SL_{3}(\bbR),\quad
\Gamma_{\infty} = \SL_{3}(\bbZ),\quad
K_{\infty} = \SO_{3}(\R),
\]
and
\begin{equation}
\label{eq:locally-symmetric-space}
\scrY = \Gamma_{\infty} \backslash G_{\infty},\qquad
\scrX = \Gamma_{\infty} \backslash G_{\infty} / K_{\infty} = \scrY/K_{\infty}.
\end{equation}
The space \(\scrX\) is an arithmetic locally symmetric space defined by the arithmetic group \(\Gamma_{\infty}\) and the symmetric space \(G_{\infty}/K_{\infty}\).
By the theorem of Borel and Harish-Chandra \cite{BorelHarishChandra1962,PlatonovRapinchuk1994}, \(\Gamma_{\infty}\) defines a unique \(G_{\infty}\)-invariant probability measure on \(\scrY\), which we denote by \(m_{\scrY}\), and let \(m_{\scrX}\) be the pushforward to \(\scrX\).
For \(g\in \GL_{3}(\R)\), we define the \emph{normalization} of \(g\) as
\[
g^{\circ} = (\det g)^{-1/3} g \in \SL_{3}(\R),
\]
using the real cube root.

Let \((M,q)\) be a positive definite ternary quadratic lattice, and let \(A\) be its Gram matrix in some \(\Z\)-basis of \(M\), i.e., the matrix of the symmetric bilinear form \(\beta\) such that \(q(x)=\beta(x,x)\).
We define the \emph{normalized Gram matrix} as \(A^{\circ}=(\det A)^{-1/3}A\), the Gram matrix of the rescaled form \((\det A)^{-1/3}q\), of determinant \(1\).
We define the \emph{shape} of the quadratic lattice \((M,q)\) as
\begin{equation}
\rmop{sh}(M,q)=\Gamma_{\infty}(A^{\circ})^{-1/2}K_{\infty} \in \scrX,
\end{equation}
where the positive square root \((A^{\circ})^{1/2}\) is defined by the polar decomposition.
Under the identification \(gK_{\infty} \mapsto (g^{-1})^{T}(g^{-1})\) of the symmetric space \(G_{\infty}/K_{\infty}\) with the set of positive symmetric matrices, the point \(\rmop{sh}(M,q)\) corresponds to the \(\Gamma_{\infty}\)-orbit \(\{\gamma^{T}A^{\circ}\gamma\mid\gamma\in \Gamma_{\infty}\}\) of \(A^{\circ}\), so is independent of the basis chosen or the scaling factor for \(q\).

We consider two families of shapes (for the zero trace and the Gross lattice).
For \(M\in \{L,\Lambda\}\), write \(M_{E}=L_{E}\) or \(M_{E}=\Lambda_{E}\) respectively.
We define the uniform and mass-weighted shape measures by
\begin{align}
\label{eq:uniform}
\nu_{p,M}
& = \frac{1}{h_{p}} \sum_{[E]\in \rmSS_{p}} \delta_{\rmop{sh}(M_{E},q_{E})},\\
\label{eq:mass}
\nu_{p,M}^{\mass}
& = \frac{1}{M_{p}} \sum_{[E]\in \rmSS_{p}} \frac{\delta_{\rmop{sh}(M_{E},q_{E})}}{|\Aut(E)|},\quad
M_{p}  = \sum_{[E]\in \rmSS_{p}} \frac{1}{|\Aut(E)|}.
\end{align}
These are probability measures on \(\scrX\), and they correspond naturally to probability measures on \(\rmSS_{p}\).
The goal of this paper is to establish (effective) equidistribution results for these probability measures, as \(p\to\infty\), in an appropriate sense.
Understanding the distribution of certain explicit families of lattices is of crucial importance in cryptography.
For these lattices to be used effectively in cryptography, they need to behave ``randomly''.
An important family of lattices that are used in cryptography are so-called \(N\)-ary lattices, for which the ``random'' behavior was established by Goldstein--Mayer \cite{GoldsteinMayer2003} using dynamics of certain Hecke operators and known results on the theory of automorphic forms.

The probability measures we consider live naturally on the suitably defined (see below) moduli space of rank three unimodular lattices.
However, because the quaternion algebra \(B_{p,\infty}\) is ramified at \(\infty\), we cannot use simple dynamics at the infinite place, and must instead appeal to \emph{adelic} results.
A significant part of this work consists in understanding these probability measures, adelically.
Establishing equidistribution in this context is, as a result, a much harder task than establishing equidistribution of \(N\)-ary lattices, and relies on the deep results of Einsiedler, Margulis, Mohammadi, and Venkatesh \cite{EinsiedlerMargulisMohammadiVenkatesh2020}.
Although our equidistribution result is stated in the archimedean setting of effective equidistribution as in earlier work of Einsiedler, Margulis, and Venkatesh \cite{EinsiedlerMargulisVenkatesh2009}, it cannot be proved by using their method, as we just explained.
Instead, we must rely on adelic effective equidistribution.

\subsection*{Archimedean Sobolev norms}

To state our main results, which are effective, we need to fix a Sobolev norm on \(\scrY\), applied to the lift of a function on \(\scrX\).
We follow \cite[\S 3.7]{EinsiedlerMargulisVenkatesh2009}.
Consider the Lie algebra \(\frakop{sl}_{3}(\R)\) of \(G_{\infty}\).
It has a \(\Q\)-structure and we first fix an integral structure \(\frakop{sl}_{3}(\Z)\).
This is a rational \(\Ad(\Gamma_{\infty})\)-stable lattice \(\frakop{sl}_{3}(\Z)\subset \frakop{sl}_{3}(\R)\) such that \([\frakop{sl}_{3}(\Z),\frakop{sl}_{3}(\Z)]\subset \frakop{sl}_{3}(\Z)\).
Fix an \(\rmop{Ad}(K_{\infty})\)-invariant Euclidean norm \(\|\;\|\) on \(\frakop{sl}_{3}(\R)\) which is sub-multiplicative under the bracket, i.e., \(\|[Z,Z^{\prime}]\|\leq \|Z\|\, \|Z^{\prime}\|\) for \(Z,Z^{\prime}\in \frakop{sl}_{3}(\R)\), and fix an orthonormal basis \(\calB_{\infty}\).
For \(y=\Gamma_{\infty}g\in \scrY\), we define its \emph{height} (see \cite[\S 3.6]{EinsiedlerMargulisVenkatesh2009})
\begin{equation}
\label{eq:archimedean-height}
\rmop{ht}(y)=\rmop{ht}_{\infty}(y)=\rmop{ht}(\Gamma_{\infty} g ) 
= \sup \bigl\{\|\rmop{Ad}(g^{-1})Z\|^{-1}\mid Z\in \frakop{sl}_{3}(\Z)\setminus\{0\}\bigr\}.
\end{equation}
This is well-defined, independent of the choice of the representative \(g\), since \(\rmop{Ad}(\Gamma_{\infty})\) preserves \(\frakop{sl}_{3}(\Z)\).
It is continuous, positive, right \(K_{\infty}\)-invariant, so \(\rmop{ht}\) also defines a function on \(\scrX\).

For \(Z\in \frakop{sl}_{3}(\R)\), let
\[
dR_{\infty}(Z) \phi(\Gamma_{\infty}g) = \left.\frac{d}{dt}\right|_{t=0}\phi\bigl(\Gamma_{\infty}g \exp(tZ)\bigr).
\]
(See \cite[Chap. II, \S 1]{Helgason2001}.)
The map \(dR_{\infty}\) is a representation of the Lie algebra \(\frakop{sl}_{3}(\R)\) on \(C^{\infty}(\scrY)\), the differential of the right regular representation of \(G_{\infty}\) on \(C^{\infty}(\scrY)\);
it extends uniquely to the universal enveloping algebra \(\frakU(\frakop{sl}_{3}(\R))\) acting on smooth functions \(C^{\infty}(\scrY)\) (\cite[Theorem 3.2.3]{Varadarajan1974}).
The universal enveloping algebra is spanned by monomials in the chosen basis.
Let \(\calD_{d}\subset \frakU(\frakop{sl}_{3}(\R))\) denote the finite set of monomials of degree at most \(d\), and write \(D\phi\) for \(dR_{\infty}(D)\phi\).
For a function \(f\) on \(\scrX\), we define the lift \(\hat{f}(y)=f(yK_{\infty})\) to \(\scrY\).
We define \(C_{c}^{\infty}(\scrX)\) to be the space of compactly supported functions \(f\) on \(\scrX\) whose lift \(\hat{f}\) is smooth on \(\scrY\).
From the definition of the Sobolev norm on \(C_{c}^{\infty}(\scrY)\), restricted to those lifts from functions on \(\scrX\), we define the \(L^{2}\)-Sobolev norms \(\calN_{d}\) on \(C_{c}^{\infty}(\scrX)\) by
\begin{equation}
\label{eq:real-Sobolev-norm}
\calN_{d}(f)^{2} = \sum_{D\in \calD_{d}} \| \rmop{ht}^{d} D \hat{f}\|_{L^{2}(\scrY,m_{\scrY})}^{2}.
\end{equation}

Note that in \cite{EinsiedlerMargulisVenkatesh2009}, the Sobolev norm is denoted by \(\calS_{d}(f)\);
we reserve that notation for the \emph{adelic} Sobolev norm to be introduced in \S \ref{subsection:adelic-Sobolev-norm}, equation \eqref{eq:adelic-Sobolev-norm}.

\subsection*{Main results}

We can now state the main results of the paper.
Our main result is an effective equidistribution result for the uniform and mass-weighted shape measures with weak convergence to the unique invariant probability measure \(m_{\scrX}\) as \(p\to\infty\) through the primes.

\begin{thm}
\label{theorem:A}
There exist an integer \(d\geq 1\), and constants \(C,\kappa>0\), such that for every prime \(p\geq 7\), every \(M\in \{L,\Lambda\}\), every \(\nu\in \{\nu_{p,M},\nu_{p,M}^{\mass}\}\), and every \(f\in C_{c}^{\infty}(\scrX)\), we have
\begin{equation}
\left|\int_{\scrX} f \, d\nu - \int_{\scrX} f\, dm_{\scrX}\right| \leq C p^{-\kappa} \calN_{d}(f).
\end{equation}
Further, all four families of probability measures in \eqref{eq:uniform} and \eqref{eq:mass} converge weakly to \(m_{\scrX}\) as \(p\to\infty\) through the primes.
\end{thm}

Once this result is established, we easily obtain several results on the limiting distribution of successive minima of the lattices \(L_{E}\) and \(\Lambda_{E}\).

For a positive definite quadratic lattice \((M,q)\) of rank \(n\), the successive minima are defined as
\begin{equation}
\label{eq:successive-minima-definition}
\lambda_{i}(M,q)=\min \bigl\{r>0 \mid \dim \rmop{span}_{\R} \{v\in M\mid q(v)\leq r^{2}\}\geq i\bigr\},\quad (1\leq i\leq n).
\end{equation}
Let \(\eta\) be the law of the successive minima of a random unimodular rank-three lattice, i.e., the pushforward of the probability measure \(m_{\scrX}\) under the map \(\scrX \to \R^{3}\), \(\Gamma_{\infty}gK_{\infty} \mapsto (\lambda_{i}(g^{-1}\Z^{3},q_{0}))_{1\leq i\leq 3}\).
By Proposition \ref{proposition:zero-trace-Gross-lattices} below, \(\det(L_{E},q_{E})=p^{2}/4\) and \(\det(\Lambda_{E},q_{E})=4p^{2}\).
We define the \emph{normalized minima} for \(M\in \{L,\Lambda\}\) as
\begin{equation}
u_{p,E}^{(M)} = \left(\frac{\lambda_{i}(M_{E},q_{E})}{\det(M_{E},q_{E})^{1/6}}\right)_{1\leq i\leq 3}.
\end{equation}
Thus,
\[
u_{p,E}^{(L)} = \left(\frac{\lambda_{i}(L_{E},q_{E})}{(p/2)^{1/3}}\right)_{1\leq i\leq 3},\quad
u_{p,E}^{(\Lambda)} = \left(\frac{\lambda_{i}(\Lambda_{E},q_{E})}{(2p)^{1/3}}\right)_{1\leq i\leq 3}.
\]
For \(M\in \{L,\Lambda\}\), we write \(u_{p,E}^{(M)}(i)\) for the \(i\)-th component of \(u_{p,E}^{(M)}\), and define
\[
\eta_{p}^{(M)} = \frac{1}{h_{p}}\sum_{[E]\in \rmSS_{p}} \delta_{u_{p,E}^{(M)}},\quad
\eta_{p}^{(M,\mass)} = \frac{1}{M_{p}}\sum_{[E]\in \rmSS_{p}} \frac{1}{|\Aut(E)|}\delta_{u_{p,E}^{(M)}},
\]
for the uniform and mass-weighted empirical laws, respectively, with \(M_{p}\) defined as in \eqref{eq:mass}.

Our next result is a precise result on the distribution of the successive minima, obtained as an easy corollary of Theorem \ref{theorem:A}.

\begin{thm}
\label{theorem:B}
\begin{enumerate}[auto]
\item For either choice of \(M\in\{L,\Lambda\}\), the empirical probability measures \(\eta_{p}^{(M)}\) and \(\eta_{p}^{(M,\mass)}\) both converge weakly to \(\eta\) as \(p\to\infty\).
Equivalently, for every Borel set \(\calR\subset\R_{>0}^{3}\) with \(\eta(\partial\calR)=0\), we have
\[
\eta_{p}^{(M)}(\calR)= \frac{\#\{[E]\in \rmSS_{p} \mid u_{p,E}^{(M)}\in \calR\}}{h_{p}} \xto[p\to\infty]{} \eta(\calR),\quad
\eta_{p}^{(M,\mass)}(\calR)\xto[p\to\infty]{} \eta(\calR).
\]
\item As \(p\to\infty\), for every function \(\psi(p)>0\) such that \(\psi(p)\to\infty\), the proportion of classes of curves satisfying, simultaneously for \(M\in\{L,\Lambda\}\), and \(i=1,2,3\),
\[
\frac{1}{\psi(p)}\leq u_{p,E}^{(M)}(i) \leq \psi(p),\quad
\frac{u_{p,E}^{(M)}(3)}{u_{p,E}^{(M)}(1)}\leq \psi(p),
\]
tends to \(1\).
This also holds for normalized mass-weighted sampling.
\end{enumerate}
\end{thm}

\section{Quaternionic lattices and groups}
\label{section:quaternion-groups}

In this section, we fix our notation and the basic objects that will be used throughout the paper.
We denote by \(\bbA\) the ring of adeles of \(\Q\), by \(\bbA_{\rmf}\) the finite adeles, and by \(\hat{\Z}\) the profinite completion of \(\Z\).
We use standard algebraic group and arithmetic group notation;
see \cite{Raghunathan1972,Margulis1991,PlatonovRapinchuk1994,KnusMerkurjevRostTignol1998} for instance.

For a group \(G\) and subgroups \(\Gamma,K\subset G\), and \(g\in G\), we generally write
\[
[g]=\Gamma gK\in \Gamma\backslash G/K,\quad
\Gamma_{g,K}=\Gamma \cap gKg^{-1}.
\]
This is the stabilizer of \(gK\) in \(\Gamma\) for the \(\Gamma\)-action on \(G/K\).

\subsection{From lattices to adelic orbits}
\label{subsection:lattices-adelic-orbits}

Let \(p\)  be a prime, and \(B=B_{p,\infty}\).
Let \(\calO\subset B\) be a maximal order.
For \(x\in B\), let \(\bar{x}=\trd(x)-x\) be quaternionic conjugation.
Let
\[
V=B^{0}=\ker(\trd:B\to\Q),\quad
q=\nrd|_{V},\quad
L=\calO\cap V,\quad
\beta(x,y)=\tfrac{1}{2}\trd(x\bar{y}).
\]
We now fix a curve \(E_{*}\in \rmSS_{p}\), and let \(\calO=\End(E_{*})\), \(B=\calO\otimes \Q\).
We then consider the algebraic \(\Q\)-group scheme
\[
\bfH=\bfO^{+}(V,q)=\bfSO(V,q).
\]

For any prime \(\ell\), we write \(B_{\ell}=B\otimes \Q_{\ell}\), \(V_{\ell} = V\otimes \Q_{\ell}\), and \(\hat{N}=\prod_{\ell} (N\otimes \Z_{\ell})\) for a full lattice \(N\) in a \(\Q\)-vector space.
In particular, this applies to orders.

We use the standard local--global dictionary for lattices in the form of adelic orbits.
Let \(V\) be a rational vector space, and \(M\subset V\) be a full lattice.
For \(h\in \bfH(\bbA_{\rmf})\), let
\begin{equation}
\label{eq:local-global-lattice}
h\cdot M = V \cap (h\hat{M}),\quad
(h\cdot M)_{\ell} = h_{\ell}M_{\ell}.
\end{equation}
The stabilizer of \(M\) in \(\bfH(\bbA_{\rmf})\) is
\[
C_{M} = \{h\in \bfH(\bbA_{\rmf}) \mid h\hat{M}=\hat{M}\}.
\]
This is a compact open subgroup.
This can be easily seen by choosing a \(\Z\)-basis of \(M\).
The faithful representation \(\bfH\into \bfGL(V)\) identifies \(C_{M}=\bfH(\bbA_{\rmf})\cap \bfGL_{3}(\hat{\Z})\).
Since \(\bfH(\bbA_{\rmf})\) is closed in \(\bfGL_{3}(\bbA_{\rmf})\) and \(\bfGL_{3}(\hat{\Z})\) is compact open, the stabilizer is compact open.
The map \(h\mapsto h\cdot M\) identifies \(\bfH(\bbA_{\rmf})/C_{M}\) with the orbit \(\bfH(\bbA_{\rmf})\cdot \hat{M}\).

\subsection{Quaternionic conjugation and norm one elements}

For a commutative \(\Q\)-algebra \(R\), we let \(B_{R}\) and \(V_{R}\) denote the extension of scalars to \(R\).
We consider the algebraic groups \(\bfGL_{1}(B)\) and \(\bfS=\bfSL_{1}(B)\) defined by their functor of points:
\[
\bfGL_{1}(B)(R)=B_{R}^{\times},\quad
\bfS(R)=\{b\in B_{R}^{\times} \mid \nrd(b)=1\}.
\]
The functor of points of \(\bfH\) is
\[
\bfH(R) = \{ h\in \SL(V_{R}) \mid q(hx)=q(x) ,\, \forall x\in V_{R}\}.
\]
Thus, \(\bfGL_{1}(B)(\Q)=B^{\times}\), \(\bfGL_{1}(B)(\Q_{\ell})=B_{\ell}^{\times}\), and \(\bfGL_{1}(B)(\bbA_{\rmf})=B(\bbA_{\rmf})^{\times}\).
We also use the projective group \(\bfPGL_{1}(B)=\bfGL_{1}(B)/\bbG_{\rmm}\) for the quotient of \(\bfGL_{1}(B)\) by its central subgroup \(\bbG_{\rmm}\).
The groups \(\bfGL_{1}(B)\), \(\bfS\), and \(\bfPGL_{1}(B)\) are the groups \(G\), \(G_{1}\), and \(G^{\mathrm{ad}}\) in \cite[\S 3, (3.1)--(3.2)]{Yu2015}.

Let \(\Psi:\bfGL_{1}(B)\to \bfH\) be the conjugation homomorphism, \(\Psi_{R}(b)(x)=bxb^{-1}\) (for \(x\in V_{R}\)).
It induces an isomorphism \(\bfPGL_{1}(B)\cong \bfH\);
this is an instance of exceptional isomorphism of groups of type \(A_{1}=B_{1}\) (see \cite[\S 26B]{KnusMerkurjevRostTignol1998}).
Let \(\rho=\Psi|_{\bfS}\).
We obtain exact sequences
\begin{equation}
\label{eq:A1-B1-exact-sequences}
1 \to \bbG_{\rmm} \to \bfGL_{1}(B) \xto{\Psi} \bfH \to 1,\quad
1 \to \bfmu_{2} \to \bfS \xto{\rho} \bfH \to 1.
\end{equation}
The first is a central extension and the second a central isogeny.
The group \(\bfS\) is an inner form of \(\bfSL_{2}\), it is connected, simply-connected, and semisimple.
Further, for any field extension \(F/\Q\), Hilbert's theorem 90 applied to the first sequence implies
\[
\bfH(F)=\Psi_{F}(B_{F}^{\times})=B_{F}^{\times} / F^{\times}.
\]
If \(c\) is a scalar and \(b\) a quaternion, since \(\nrd(cb)=c^{2}\nrd(b)\), we see that an element \(\Psi_{F}(b)\) belongs to \(\rho(\bfS(F))\) exactly when \(\nrd(b)\in F^{\times 2}\).
\[
\rho(\bfS(F))=\{\Psi_{F}(b) \mid b\in B_{F}^{\times},\, \nrd(b)\in F^{\times 2}\}.
\]
We have the diagram of \(\Q\)-group schemes:
\[
\begin{tikzcd}
1 \rar & \bfmu_{2} \dar \rar &  \bfS\rar{\rho} \dar[hookrightarrow]& \bfH\dar[equal] \rar &1\\
1 \rar & \bbG_{\rmm} \rar & \bfGL_{1}(B)\dar[twoheadrightarrow] \rar{\Psi} & \bfH\dar[equal]\rar & 1\\
&& \bfPGL_{1}(B) \rar{\cong} & \bfH
\end{tikzcd}
\]

\subsection{Orders and integral levels}

We have defined the \(\Q\)-groups \(\bfGL_{1}(B)\) and \(\bfS\).
We use the following integral structure.
We define affine \(\Z\)-group schemes
\[
\calop{GL}_{1}(\calO)(R)=(\calO \otimes_{\Z}R)^{\times},\quad
\calS(R) = \{s\in \calop{GL}_{1}(\calO)(R)\mid s\bar{s}=1\},
\]
for commutative rings \(R\).
Their generic fibers are \(\bfGL_{1}(B)\) and \(\bfS\) over \(\Q\), respectively.
In \cite[\S 2.4]{Yu2005}, the group scheme \(G_{x}\) is defined from the Rosati involution;
for the canonical principal polarization on supersingular elliptic curves, the Rosati involution coincides with quaternion conjugation, so the group scheme \(\calS\) we defined is the group scheme \(G_{x}\) from \cite[\S 2.4]{Yu2005}.
We define the integral levels as in \cite[p.~167 and \S 3.3, p.~177]{Yu2015}.
We define
\[
U^{1} = \calS(\hat{\Z})= \prod_{\ell} \calO_{\ell}^{1},\quad
\calO_{\ell}^{1} = \{s\in \calO_{\ell}^{\times} \mid \nrd(s)=1\}.
\]
The group \(U^{1}\) is a compact open subgroup of \(\bfS(\bbA_{\rmf})\);
it is denoted by \(U_{1}=\hat{\calO}{}^{\times}\cap G_{1}(\bbA_{\rmf})\) in \cite{Yu2015}.
We use standard facts about the local description of maximal quaternion orders (see \cite{Reiner1975,MaclachlanReid2003,Vigneras1980}).

For \(s\in \bfS(\bbA_{\rmf})\), define
\begin{equation}
\label{eq:conjugate-stabilizer-order}
\calO_{s} = B \cap s\hat{\calO} s^{-1},\quad
(\calO_{s})_{\ell} = s_{\ell} \calO_{\ell} s_{\ell}^{-1}.
\end{equation}
This is a well-defined order.

\subsection{The Gross lattice}

Let \(\calO\) be a maximal order in \(B\), and let \(V=B^{0}\).
We write in this section \(T(x)=2x-\trd(x)\) for \(x\in \calO\).
The Gross lattice is
\begin{equation}
\Lambda = T(\calO) = (\Z+2\calO)\cap V,
\end{equation}
with quadratic form \(q=\nrd|_{\Lambda}\);
see \cite[\S 3.3]{GorenLove2025} for properties of the Gross lattice.
We gather the relevant facts in the following proposition.

\begin{proposition}
\label{proposition:zero-trace-Gross-lattices}
The lattices \(L\) and \(\Lambda\) satisfy the following properties:
\begin{alignat}{2}
2L \subset \Lambda & \subset L,&\quad
[L:\Lambda] & = 4,\quad
[\Lambda:2L]=2,\\
\det(L,q) & = p^{2}/4,&\quad
\det(\Lambda,q)&=4p^{2}.
\end{alignat}
\end{proposition}

We verify the compatibility of the lattices \(L\) and \(\Lambda\) with (adelic) conjugation.
Since \(\trd(bxb^{-1})=\trd(x)\) and \(T(bxb^{-1})=bT(x)b^{-1}\), we obtain from \eqref{eq:conjugate-stabilizer-order}
\begin{align*}
(\calO_{s} \cap V)_{\ell}
& = \ker(\trd|_{(\calO_{s})_{\ell}})
 = \rho(s_{\ell})L_{\ell},\\
T(\calO_{s})_{\ell}
& = T(s_{\ell}\calO_{\ell}s_{\ell}^{-1})
 = \rho(s_{\ell}) \Lambda_{\ell}.
\end{align*}
From \eqref{eq:local-global-lattice}, it follows that
\begin{equation}
\calO_{s} \cap V = \rho(s) \cdot L,\quad
T(\calO_{s}) = \rho(s)\cdot \Lambda.
\end{equation}

\section{Supersingular curves and the norm-one group}

Let \(\scrC_{\bfS}=\bfS(\Q)\bs \bfS(\bbA_{\rmf}) / \calS(\hat{\Z})\).
The Deuring correspondence is usually stated using the double quotient \(\bfGL_{1}(B)(\Q)\bs \bfGL_{1}(B)(\bbA_{\rmf}) / \calop{GL}_{1}(\calO)(\hat{\Z})\).
But we will recast the Deuring correspondence directly in the context of the group \(\bfS\) and the lattices \(L\) and \(\Lambda\), using Yu's parametrization (see \cite{Yu2005,Yu2015}).

\subsection{Yu's parametrization and the mass}

In this section, we state the relevant results from \cite{Yu2005} (which refers to \cite{Yu2003}).
Yu's parametrization is given for polarized abelian varieties.
Every elliptic curve \(E\) has a canonical principal polarization \(\lambda_{E}\), and in dimension \(1\), being supersingular is equivalent to being superspecial.
So for \(x=(E_{*},\lambda_{E_{*}})\), Yu's geometric class set (see \cite[\S 2.4]{Yu2005}) is precisely \(\rmSS_{p}\).

\begin{theorem}
\label{theorem:yu-norm-one-correspondence}
Let \(\scrC_{\bfS}=\bfS(\Q)\bs \bfS(\bbA_{\rmf}) / U^{1}\).
This double-coset class set is finite, and there is a bijection
\begin{equation}
\label{eq:bijection-SSp-double-quotient-norm-one}
\Theta:\rmSS_{p} \stosim \scrC_{\bfS},\quad
\Theta([E_{*}])=[1].
\end{equation}
If \(s\in \bfS(\bbA_{\rmf})\) represents \(\Theta([E])\), then we have \(\Aut(E) \cong \bfS(\Q) \cap sU^{1}s^{-1}\).
Moreover,
\begin{equation}
\label{eq:mass-formula-Aut}
M_{p} = \sum_{[E]\in \rmSS_{p}} \frac{1}{|\Aut(E)|} = \rmop{Mass}(\bfS,U^{1}) = \frac{p-1}{24}.
\end{equation}
\end{theorem}

\begin{proof}
This follows from Theorem 2.1, Corollary 2.3, Lemma 2.4, from \cite{Yu2005}.
The mass formula follows from specialization of \cite[(3.18)]{Yu2015}, where it reads
\[
\rmop{Mass}(\bfS,U^{1}) = -\tfrac{1}{2} \zeta(-1)(p-1) = \frac{p-1}{24}.\qedhere
\]
\end{proof}

\subsection{Associated orders and lattices}

We now extend our parametrization to lattices.

\begin{proposition}
\label{proposition:stabilizer-conjugate}
Let \(s\in \bfS(\bbA_{\rmf})\), and let \(E\) represent the isomorphism class \(\Theta^{-1}([s])\).
There is an isomorphism of quaternion orders
\begin{equation}
\label{eq:isomorphism-conjugation-endomorphism-ring}
\End(E) \cong \calO_{s},\quad
\calO_{s}^{\times} = \bfS(\Q) \cap sU^{1}s^{-1}.
\end{equation}
Moreover, the corresponding isometries
\[
(L_{E},q_{E}) \simeq (\rho(s) \cdot L,q),\quad
(\Lambda_{E},q_{E}) \simeq (\rho(s) \cdot \Lambda,q),
\]
are compatible with the lattice inclusions.
Their isometry classes depend only on \([s]\).
\end{proposition}

\begin{proof}
By \cite[Proposition 2.2]{Yu2005}, there is \(\varphi\in \rmop{Hom}(E_{*},E)\otimes\Q\) such that \(\varphi^{*}\lambda_{E}=\lambda_{E_{*}}\).
Such \(\varphi\) is invertible in \(\rmop{Hom}(E_{*},E)\otimes\Q\).
Thus, \(\calO^{\prime}=\varphi^{-1}\End(E)\varphi\subset B\) is a maximal order.
For every prime \(\ell\), we have \((\calO^{\prime})_{\ell}^{1} = s_{\ell}\calO_{\ell}^{1} s_{\ell}^{-1}\) from \cite[Proof of Lemma 2.4]{Yu2005}.
For unramified \(\ell\neq p\), a maximal order of \(B_{\ell}\cong \rmop{M}_{2}(\Q_{\ell})\) is conjugate to \(\rmop{M}_{2}(\Z_{\ell})\) and spanned by norm-one units, which implies local equality of the orders from the norm-one elements:
\((\calO^{\prime})_{\ell} = s_{\ell}\calO_{\ell} s_{\ell}^{-1}\).
Recovering the global order from the completions, this proves the first part of \eqref{eq:isomorphism-conjugation-endomorphism-ring}.

A unit \(b\in \calO_{s}\) must satisfy \(\nrd(b)=1\) since \(\nrd(b),\nrd(b^{-1})\in \Z_{>0}\).
It follows from \eqref{eq:conjugate-stabilizer-order} that
\[
\calO_{s}^{\times} = \{b\in \bfS(\Q) \mid b\in s_{\ell}\calO_{\ell}^{1} s_{\ell}^{-1},\, \forall \ell\} = \bfS(\Q)\cap sU^{1}s^{-1}.
\]
Conjugation by \(\varphi\) preserves \(\trd\) and \(\nrd\), so it carries \(L_{E}=\ker(\trd|_{\End(E)})\) isometrically onto \(\calO_{s} \cap V\) and \(\Lambda_{E} = T(\End(E))\) onto \(T(\calO_{s})\), which are \(\rho(s)\cdot L\) and \(\rho(s)\cdot \Lambda\) respectively, and the inclusion is respected.
Finally, if we replace \(s\) by \(s^{\prime}=\gamma s u\) with \(\gamma \in \bfS(\Q)\) and \(u\in U^{1}\), we see that \(\rho(s)\cdot M\) is changed to \(\rho(\gamma)(\rho(s)\cdot M)\), and \(\rho(\gamma)\) is a rational isometry.
So the isometry class does not depend on the representative.
\end{proof}

We conclude that the mass-weighted shape measure \eqref{eq:mass} is, for \(M\in \{L,\Lambda\}\),
\begin{equation}
\label{eq:mass-weighted-measure}
\nu_{p,M}^{\mass} = \frac{1}{M_{p}} \sum_{[s]\in \scrC_{\bfS}} \frac{\delta_{\rmop{sh}(\rho(s)\cdot M,q)}}{|\calO_{s}^{\times}|}.
\end{equation}
This transports the mass-weighted shape measure via \(\Theta\).

\section{Homogeneous dynamics}

A significant part of this work consists in constructing the correct objects that we can understand via homogeneous dynamics.
In this section all topological groups are locally compact Hausdorff.
We let \(M\in\{L,\Lambda\}\), with \(M_{E}\) as the corresponding lattice arising explicitly from the Deuring correspondence.
We will establish the following diagram.
\begin{equation}
\label{eq:goal-diagram}
\begin{tikzcd}[row sep=3em,column sep=4em]
\rmSS_{p}\dar[swap]{[E]\mapsto \rmop{sh}(M_{E},q_{E})} & \scrC_{\bfS}\lar[swap]{\Theta^{-1}} & \bfS(\Q)\bs \bfS(\bbA)\lar{(\;)/W}[swap]{\rmop{cl}}\dar{\Phi_{p}}\\
\scrX &\scrY \cong X_{\bbA}/K_{\rmf} \lar[swap]{(\;)/K_{\infty}}& X_{\bbA}\lar[swap]{(\;)/K_{\rmf}}
\end{tikzcd}
\end{equation}
Here, \(W=\bfS(\R)\times U^{1}\) and \(\rmop{cl}(\bfS(\Q)s)=[s_{\rmf}]\).

For a unimodular locally compact group \(J\) with Haar measure \(dg\) and a discrete subgroup \(\Delta\subset J\), there is a unique \(J\)-invariant measure \(d\bar{g}\) on \(\Delta \backslash J\) such that
\begin{equation}
\label{eq:Weil}
\int_{J} \psi(g) \, dg = \int_{\Delta \backslash J} \sum_{\delta\in \Delta} \psi(\delta g) d\bar{g},\quad \psi\in C_{c}(J).
\end{equation}
(See \cite[Theorem 2.51]{Folland2016}.)
Whenever the quotient has finite volume, we rescale the quotient measure to a probability measure.

\subsection{The ambient quotient}

Let
\begin{equation}
\label{eq:PGL}
\bfG=\bfPGL_{3},\quad
X_{\bbA} = \bfG(\Q)\bs \bfG(\bbA),\quad
K_{\rmf}=\prod_{\ell} \bfPGL_{3}(\Z_{\ell}).
\end{equation}
The group \(\Z\)-scheme \(\bfPGL_{n,\Z}\) is the quotient \(\bfGL_{n,\Z}/\bbG_{\rmm}\);
see \cite[II, \S 1, no. 2.6 and III, \S 3, no. 3.4]{DemazureGabriel1970}.
This implies that for every local ring \(R\), we have
\begin{equation}
\label{eq:PGLnZ}
\bfPGL_{n}(R)=\GL_{n}(R)/R^{\times}.
\end{equation}
(See \cite[Chapter IV, Corollary 1.5]{Milne1980}.)
In particular, \(K_{\rmf}\) is compact open for the restricted product topology, and we also get \(\bfG(\bbA)=\GL_{3}(\bbA)/\bbA^{\times}\).

To define the adelic Sobolev norms as in \cite{EinsiedlerMargulisMohammadiVenkatesh2020}, we need a faithful representation of \(\bfG\).
Although we cannot use the adjoint representation over \(\Z\), a natural choice is the following explicit conjugation representation
\begin{equation}
\label{eq:faithful-representation}
\tau:\bfG \into \bfSL(\frakop{sl}_{3,\Z}) \cong \bfSL_{8},\quad
\bbP(g) \mapsto \bigl(Z\mapsto gZg^{-1}\bigr).
\end{equation}
We can choose for \(\frakop{sl}_{3,\Z}\) the standard Chevalley basis (see \cite[\S 1, Theorem 1, Remarks (a) and (c)]{Steinberg2016}).
The conjugation homomorphism \(\bfGL_{3,\Z} \to \bfGL(\frakop{sl}_{3,\Z})\) has determinant one for \(\bfGL_{3,\Z}\) and with kernel \(\bbG_{\rmm,\Z}\) (for an \(R\)-point lies in the kernel if and only if it commutes with every \(X\in \frakop{sl}_{3}(R)\), hence is scalar), so it factors through the quotient \(\bfPGL_{3,\Z}\).
It descends to a representation \(\tau:\bfPGL_{3,\Z} \to \bfGL(\frakop{sl}_{3,\Z})\) with trivial kernel, hence \(\tau_{\Q}\) is a closed immersion \cite[Theorem 5.34]{Milne2017}.
Up to the canonical identification \(\frakop{sl}_{3,\Q}\stosim \frakop{pgl}_{3,\Q}\), the generic fiber \(\tau=\tau_{\Q}\) coincides with the adjoint representation of \(\bfG\) (the differential of projectivization intertwines the representations: \(d\bbP \circ \tau(g) = \rmop{Ad}_{\bfG}(g)\circ d\bbP\)).

\begin{lemma}
\label{lemma:class-number-one}
We have \(\bfG(\bbA_{\rmf})=\bfG(\Q)K_{\rmf}\) and \(\bfG(\Q)\cap K_{\rmf}=\bbP(\GL_{3}(\Z))\).
\end{lemma}

\begin{proof}
The first assertion \(\bfG(\bbA_{\rmf})=\bfG(\Q)K_{\rmf}\) means precisely that \(\bfG\) has class number one.
By \cite[Proposition 2.2]{Borel1963}, the class number of \(\bfGL_{n,k}\) is the class number of \(k\), which is one for \(k=\Q\).
Thus, we have \(\bfGL_{3}(\bbA_{\rmf})=\bfGL_{3}(\Q) \bfGL_{3}(\hat{\Z})\), and projectivization implies \(\bfG(\bbA_{\rmf})=\bfG(\Q)K_{\rmf}\), since the maps \(\bfGL_{3}(\Q_{\ell}) \to \bfG(\Q_{\ell})\) and \(\bfGL_{3}(\Z_{\ell}) \to \bfPGL_{3}(\Z_{\ell})\) are surjective.

For the second assertion, we have by \eqref{eq:PGLnZ} that \(\bfG(\Q)=\bbP(\GL_{3}(\Q))\).
Let \(x\in \GL_{3}(\Z)\).
Then, \(\bbP(x)\in \bfG(\Q)\), but also \(x\in \GL_{3}(\Z_{\ell})\) for every prime \(\ell\), so \(\bbP(x)\in \bfPGL_{3}(\Z_{\ell})\) for every \(\ell\), hence \(\bbP(x)\in K_{\rmf}\).
This proves the inclusion \(\supset\).
We now prove the reverse inclusion \(\subset\).
Let \(x\in \bfG(\Q)\cap K_{\rmf}\).
From \eqref{eq:PGLnZ}, there is \(g\in \GL_{3}(\Q)\) such that \(x=\bbP(g)\).
We can write \(g=\varepsilon k\) where \(\varepsilon\in \bbA_{\rmf}^{\times}\) and \(k\in \GL_{3}(\hat{\Z})\).
At the primes \(\ell\) where \(g\in \GL_{3}(\Z_{\ell})\), pick \(\varepsilon_{\ell}=1\) and \(k_{\ell}=g\).
This holds for almost all primes \(\ell\), hence \(\varepsilon=(\varepsilon_{\ell})_{\ell}\in \bbA_{\rmf}^{\times}\) and \(k=(k_{\ell})_{\ell}\in \GL_{3}(\hat{\Z})\).
Since \(\bfGL_{1,\Q}\) has class number one, we can write \(\varepsilon=cu\) with \(c\in \Q^{\times}\) and \(u\in \hat{\Z}{}^{\times}\).
Then, \(c^{-1}g=uk\in \GL_{3}(\Q)\cap \GL_{3}(\hat{\Z})=\GL_{3}(\Z)\).
Thus, \(x\in \bbP(\GL_{3}(\Z))\).
\end{proof}

We now choose once and for all a \(\Z\)-basis of the lattice \(M\), identifying \(V\) with \(\Q^{3}\), and \(M\) with \(\Z^{3}\).
Let \(\Xi=\Xi_{p,M}\) be the Gram matrix of the bilinear form \(\beta\) defined in \S \ref{subsection:lattices-adelic-orbits}, in this basis, so the matrix \(\Xi\) is \((\beta(v_{i},v_{j}))_{i,j}\).
This gives a fixed coordinate realization of \(\bfH=\bfSO(V,q)\):
\[
\bfH(R)=\{h\in \SL_{3}(R) \mid h^{T} \Xi h=\Xi\},
\]
for every commutative \(\Q\)-algebra \(R\).
This makes \(\bfH\) a closed subgroup scheme of \(\bfSL_{3}\).
This coordinate realization is crucial because the ambient data must be independent of \(p\), whereas \(V\), \(q\), \(L\) are not.

Recall the map \(\rho:\bfS \to \bfH\).
We now define
\begin{equation}
\label{eq:future-datum}
\left\{\begin{alignedat}{2}
\iota=\iota_{M}:\bfH&\hookrightarrow \bfSL_{3} \to \bfG,&\quad
j=j_{M}&=\iota\circ \rho:\bfS \to \bfG,\\
\xi =\xi_{p,M}& = (\Xi^{\circ})^{-1/2}\in \SL_{3}(\R),&\quad
\omega_{p}=\omega_{p,M}&=(\bbP(\xi),1_{\rmf})\in \bfG(\bbA).
\end{alignedat}\right.
\end{equation}

Since \(\ker(\bfGL_{3} \to \bfPGL_{3})=\bbG_{\rmm}\) via scalar matrices, \(\ker\iota=\bfH\cap \bbG_{\rmm}\).
Consider a scalar matrix \(zI_{3}\in \bfH(R)\);
then \(z^{2}\Xi=\Xi\), hence \(z^{2}=1\).
The determinant one condition implies \(z^{3}=1\), and together we get \(z=1\), hence \(\ker\iota=\{1\}\).
Thus, \(\iota\) is a monomorphism, hence a closed immersion of algebraic groups \cite[Proposition 5.31, Theorem 5.34]{Milne2017}.
This identifies \(\bfH\)  with a closed subgroup scheme \(\iota(\bfH)\) of \(\bfG\) defined over \(\Q\), and since \(\ker\iota=1\), we have \(\ker j = \ker \rho = \bfmu_{2}\), as \(\Q\)-group schemes.
Moreover, since \(\sigma=(\Xi^{\circ})^{1/2}\) is symmetric and satisfies \(\sigma^{2}=\Xi^{\circ}\), we have, for \( \xi =\sigma^{-1}\), \( \xi ^{T} \Xi^{\circ}  \xi =\sigma^{-1}\sigma^{2}\sigma^{-1}=I_{3}\).
Since we can rewrite \(\bfH(\R)=\{h\in \SL_{3}(\R) \mid h^{T}\Xi^{\circ}h=\Xi^{\circ}\}\), we can also rewrite \(\bfH(\R)\) in terms of \( \xi \):
we have \(\Xi^{\circ} = ( \xi ^{-1})^{T}( \xi ^{-1})\), so
\begin{align*}
h^{T}\Xi^{\circ}h=\Xi^{\circ}
&\quad\Leftrightarrow\quad
h^{T}( \xi ^{-1})^{T}( \xi ^{-1})h=( \xi ^{-1})^{T}( \xi ^{-1})\\
&\quad\Leftrightarrow\quad
( \xi ^{-1}h \xi )^{T}( \xi ^{-1}h \xi )=I_{3}\\
&\quad\Leftrightarrow\quad
( \xi ^{-1}h \xi )\in \rmop{O}_{3}(\R).
\end{align*}
Since \(\det(h)=1\) if and only if \(\det( \xi ^{-1}h \xi )=1\), we have
\[
h\in \bfH(\R)
\quad\Leftrightarrow\quad
( \xi ^{-1}h \xi )\in \rmop{SO}_{3}(\R)=K_{\infty}.
\]
Thus, \( \xi ^{-1} \bfH(\R) \xi =K_{\infty}\).
We summarize:
\begin{equation}
\label{eq:kernels-compact-subgroup}
\ker\iota=1,\quad
\ker j = \bfmu_{2},\quad
 \xi ^{-1} \bfH(\R) \xi  = K_{\infty}.
\end{equation}

We have an archimedean quotient, but the equidistribution theorem is adelic.
The next proposition gives the adelic realization of the archimedean shape space.

Let \(\bfG(\R)\) act on \(X_{\bbA}\) on the right by \(h\mapsto (h,1_{\rmf})\);
this action commutes with that of \(K_{\rmf}\) so factors through the quotient to \(X_{\bbA}/K_{\rmf}\).
We can identify \(\bfG(\R)\) with \(G_{\infty}\) by the isomorphism of Lie groups \(\bbP(g) \mapsto g^{\circ}\).

\begin{proposition}
\label{proposition:adelic-shape-space}
Let \(\Delta_{\infty} = \bfG(\Q) \cap K_{\rmf} = \bbP(\GL_{3}(\Z))\).
For \(a=(a_{\infty},a_{\rmf})\in \GL_{3}(\bbA)\), let \(N(a_{\rmf})=a_{\rmf} \cdot \Z^{3} = \Q^{3}\cap a_{\rmf}\hat{\Z}{}^{3}\), and define the quadratic forms \(q_{0}\) and \(q_{a_{\infty}}\) on \(\R^{3}\) by
\[
q_{0}(v)=v^{T}v,\quad
q_{a_{\infty}}(v) = q_{0}(a_{\infty}^{-1}v),\quad v\in \R^{3}.
\]
\begin{enumerate}[auto]
\item 
The maps
\begin{alignat*}{2}
\Delta_{\infty} \bs (\bfG(\R) \times K_{\rmf}) & \to X_{\bbA},&\quad
\Delta_{\infty}(h,k_{\rmf}) &\mapsto \bfG(\Q) (h,k_{\rmf}),\\
\Delta_{\infty}\bs \bfG(\R) & \to X_{\bbA}/K_{\rmf},&\quad
\Delta_{\infty} h& \mapsto \bfG(\Q)(h,1_{\rmf}) K_{\rmf},
\end{alignat*}
are \(\bfG(\R)\)-equivariant homeomorphisms.
\item There are right \(\bfG(\R)\)-equivariant homeomorphisms
\begin{equation}
\label{eq:adelic-real-isomorphism}
X_{\bbA} /K_{\rmf} \stosim \Delta_{\infty} \bs \bfG(\R) \stosim \scrY,
\end{equation}
given by
\[
\bfG(\Q)(h,1_{\rmf})K_{\rmf} \mapsto \Delta_{\infty} h,\quad
\Delta_{\infty} \bbP(g) \mapsto \Gamma_{\infty} g^{\circ},\quad
(h\in \bfG(\R), g\in \GL_{3}(\R)).
\]
Composing this with the quotient map \(X_{\bbA} \to X_{\bbA}/K_{\rmf}\) defines \(r:X_{\bbA} \to \scrY\), satisfying the following.
For \(x=\bfG(\Q)(\bbP(a),k_{\rmf})\in X_{\bbA}\) with \(k_{\rmf}\in K_{\rmf}\) and \( a\in \GL_{3}(\R)\), and for any \(u\in K_{\rmf}\), we have \(r(xu)=r(x)=\Gamma_{\infty}a^{\circ}\), and \(r(x(\bbP(a^{\prime}),1_{\rmf})) = r(x)a^{\prime\circ}\) for any \(a^{\prime}\in \GL_{3}(\R)\).

\item The maps \(r:X_{\bbA} \to \scrY\) and \(\rmop{pr}:X_{\bbA} \to \scrX\), \(\pr(x)=r(x)K_{\infty}\), are continuous and proper.

\item \(X_{\bbA}\) carries a unique \(\bfG(\bbA)\)-invariant probability measure \(m_{\bbA}\), and
\begin{equation}
\label{eq:pushforward-adelic-measure}
r_{*} m_{\bbA} = m_{\scrY},\quad
\rmop{pr}_{*} m_{\bbA}=m_{\scrX}.
\end{equation}

\item 
We restrict \(q_{a_{\infty}}\) to the lattice \(N(a_{\rmf})\subset \Q^{3}\subset \R^{3}\).
Then, for every \(a\in \GL_{3}(\bbA)\),
\begin{equation}
\label{equation:projection-shape-lattice-correspondence}
\rmop{pr}\bigl(\bfG(\Q)\bbP(a)\bigr)
=\rmop{sh}\bigl(N(a_{\rmf}),q_{a_{\infty}}\bigr)
=\rmop{sh}\bigl(a_{\infty}^{-1}N(a_{\rmf}),q_{0}\bigr).
\end{equation}
\item For \(f\in C_{c}^{\infty}(\scrX)\), the function \(f\circ \rmop{pr}\) is smooth, compactly supported, and right \(K_{\rmf}\)-invariant on \(X_{\bbA}\), that is, \(f\circ \rmop{pr}\in C_{c}^{\infty}(X_{\bbA})^{K_{\rmf}}\), and \(f\circ \pr = \hat{f}\circ r\).
\end{enumerate}
\end{proposition}

\begin{proof}
(1)
The first map is continuous and \(\bfG(\R)\)-equivariant by construction.
By Lemma \ref{lemma:class-number-one}, \(\bfG(\bbA) = \bfG(\Q)(\bfG(\R)\times K_{\rmf})\), so it is surjective, and \(\bfG(\Q)\cap (\bfG(\R)\times K_{\rmf})=\Delta_{\infty}\), diagonally embedded, so it is injective.
It is an open map since \(\bfG(\R)\times K_{\rmf}\) is open in \(\bfG(\bbA)\) and the quotient map \(\bfG(\bbA)\to X_{\bbA}\) is open.
Quotienting by the right action of \(K_{\rmf}\)  gives the second homeomorphism, equivariant by construction.

(2) In odd dimension, the map \(\bbP(g) \mapsto g^{\circ}\) is an isomorphism of Lie groups \(\bfG(\R)\stosim G_{\infty}\).
Since \(\Delta_{\infty} = \bbP(\GL_{3}(\Z))\), this maps \(\Delta_{\infty}\) onto \(\Gamma_{\infty}\) (since \(\gamma^{\circ} = (\det\gamma)\gamma\) for \(\gamma\in\GL_{3}(\Z)\)), so the map restricts to an isomorphism \(\Delta_{\infty} \stosim \Gamma_{\infty}\), and it induces the second map \(\Delta_{\infty} \bs \bfG(\R) \stosim \scrY\).

For the formula, for \(x\in X_{\bbA}\), we can write \(x=\bfG(\Q)(h,k)\) with \(h\in \bfG(\R)\) and \(k\in K_{\rmf}\).
The right action of \(\bfG(\bbA)\) on \(X_{\bbA}\), restricted to \(\bfG(\R)\), gives
\[
x(\bbP(a),1_{\rmf})=\bfG(\Q)(h\bbP(a),k)=\bfG(\Q)(\bbP(ga),k),
\]
where we can choose \(g\in \GL_{3}(\R)\) with \(h=\bbP(g)\).
Then, \(r(x)=\Gamma_{\infty}g^{\circ}\in \scrY\), so
\[
r\bigl(x(\bbP(a),1_{\rmf})\bigr)=r\bigl(\bfG(\Q)(\bbP(ga),k)\bigr)=\Gamma_{\infty}(ga)^{\circ}
=\Gamma_{\infty}g^{\circ}a^{\circ}=r(x)a^{\circ}.
\]

(3) The quotient map by a compact group is continuous and proper;
the map \(r\) is the quotient map by \(K_{\rmf}\) followed by the homeomorphism \eqref{eq:adelic-real-isomorphism}, and \(\pr\) is \(r\) followed by the quotient by \(K_{\infty}\).

(4) The uniqueness of \(m_{\bbA}\) is standard.
By the equivariance in (2), \(r_{*}m_{\bbA}\) is a \(G_{\infty}\)-invariant probability measure on \(\scrY\), so it equals \(m_{\scrY}\) by uniqueness, and since \(m_{\scrX}\) is the pushforward of \(m_{\scrY}\) under the quotient map, the pushforward \(\pr_{*}m_{\bbA}=m_{\scrX}\) follows.

(5) First assume that \(a_{\rmf}\in \GL_{3}(\hat{\Z})\) is integral, so write \(a=(g,k_{\rmf})\) with \(g\in \GL_{3}(\R)\) and \(k_{\rmf}\in \GL_{3}(\hat{\Z})\).
Then, \(a_{\rmf}\hat{\Z}{}^{3}=\hat{\Z}{}^{3}\), so \(N(a_{\rmf})=\Z^{3}\).
Put \(u=g^{\circ}\).
Write \(x=\bfG(\Q)\bbP(a)\).
Since \(\bbP(k_{\rmf})\in K_{\rmf}\), we can write \(x=x_{0}\cdot (1,\bbP(k_{\rmf}))\) with \(x_{0}=\bfG(\Q)(\bbP(g),1_{\rmf})\).
Since \(r\) is right \(K_{\rmf}\)-invariant, \(r(x)=r(x_{0})=\Gamma_{\infty}g^{\circ}=\Gamma_{\infty} u\) by \eqref{eq:adelic-real-isomorphism}, hence \(\pr(x)=\Gamma_{\infty} u K_{\infty}\).
Now, recall that for a quadratic lattice \((M,q)\), \(\sh(M,q)=\Gamma_{\infty}( A ^{\circ})^{-1/2}K_{\infty}\) where \( A \) is the Gram matrix of the form \(q\) in the fixed basis of \(M\).
If \((e_{i})\) is the fixed basis, then \((u^{-1}e_{i})\) is a basis of \(u^{-1}\Z^{3}\).
Since \(q_{0}=v^{T}v\), the Gram matrix is then \( A _{ij}=(u^{-1}e_{i})^{T}(u^{-1}e_{j})=e_{i}^{T}(u^{-T}u^{-1})e_{j}\), so that \( A =u^{-T}u^{-1}=(uu^{T})^{-1}\).
Since \(u=g^{\circ}\), we have \(\det(u)=1=\det( A )\), hence \( A ^{\circ}= A \).
It follows that
\[
\sh(u^{-1}\Z^{3},q_{0}) = \Gamma_{\infty} A ^{-1/2}K_{\infty}.
\]
Now let \(k= A ^{1/2}u\).
Since \( A \) is symmetric, \(kk^{T}= A ^{1/2}uu^{T} A ^{1/2}= A ^{1/2} A ^{-1} A ^{1/2}=I_{3}\), and \(\det k=1\), hence \(k\in K_{\infty}=\SO_{3}(\R)\).
Thus, \(u= A ^{-1/2}k\in  A ^{-1/2}K_{\infty}\) and \( A ^{-1/2}K_{\infty} =uK_{\infty}\).
The map \(v\mapsto g^{-1}v\) is an isometry \((\Z^{3},q_{g})\to (g^{-1}\Z^{3},q_{0})\), and the lattices \(u^{-1}\Z^{3}\) and \(g^{-1}\Z^{3}\) are homothetic.
We conclude that
\[
\sh(\Z^{3},q_{g}) = \sh(u^{-1}\Z^{3},q_{0})=\Gamma_{\infty} u K_{\infty} = \Gamma_{\infty} g^{\circ} K_{\infty}.
\]
Thus,
\[
\pr(\bfG(\Q)\bbP(a))=\pr(x)=\Gamma_{\infty} g^{\circ} K_{\infty} = \sh(g^{-1}\Z^{3},q_{0})=\sh(\Z^{3},q_{g}).
\]

Now for the general case, let \(a=(a_{\infty},a_{\rmf})\in \GL_{3}(\bbA)\), let \(\gamma\in\GL_{3}(\Q)\) be a basis matrix of the lattice \(N=N(a_{\rmf})\), so that \(N=\gamma \Z^{3}\).
Taking completions, \(\gamma\hat{\Z}{}^{3}=\hat{N}=a_{\rmf}\hat{\Z}{}^{3}\), with the diagonal embedding of \(\gamma\) in the finite adeles, so \(\gamma^{-1}a_{\rmf}\in \GL_{3}(\hat{\Z})\).
Let \(b_{\infty}=\gamma^{-1}a_{\infty}\).
Since \(\bbP(\gamma)\in \bfG(\Q)\), we have
\begin{align*}
\bfG(\Q)\bbP(a)
& = \bfG(\Q) (\bbP(a_{\infty}),\bbP(a_{\rmf}))\\
& = \bfG(\Q) \bbP(\gamma^{-1})(\bbP(a_{\infty}),\bbP(a_{\rmf}))
 = \bfG(\Q) (\bbP(b_{\infty}),\bbP(\gamma^{-1}a_{\rmf})),
\end{align*}
and the integral case shows that
\[
\pr(\bfG(\Q)\bbP(a))=\sh(b_{\infty}^{-1}\Z^{3},q_{0})
=\sh(a_{\infty}^{-1}\gamma\Z^{3},q_{0})
=\sh(a_{\infty}^{-1}N(a_{\rmf}),q_{0}).
\]
Since \(v\mapsto a_{\infty}^{-1}v\) is an isometry \((N,q_{a_{\infty}})\to (a_{\infty}^{-1}N,q_{0})\), this completes the proof of (5).

(6) Let \(f\in C_{c}^{\infty}(\scrX)\), i.e., we have defined \(\hat{f}(y)=f(yK_{\infty})\).
This is a smooth function on \(\scrY\), and for every \(x\in X_{\bbA}\), we have \((f\circ \pr)(x)=f(\pr(x))=f(r(x)K_{\infty})=\hat{f}(r(x))=(\hat{f}\circ r)(x)\).
Thus, \(f\circ \pr=\hat{f}\circ r\), and it is \(K_{\rmf}\)-invariant by (2).
Let \(F=f\circ \pr\).
Under the homeomorphism of (1), if we write \(\pi(h,k_{\rmf})=\bfG(\Q)(h,k_{\rmf})\), with \(h=\bbP(g)\), we get \((F\circ \pi)(\bbP(g),k_{\rmf})=\hat{f}(\Gamma_{\infty}g^{\circ})\).
It follows that \(F\) is smooth on \(X_{\bbA}\), with compact support \(\pr^{-1}(\rmop{supp}f)\), since \(\pr\) is proper.
\end{proof}

\begin{lemma}
\label{lemma:integral-structures}
For either case \(M\in \{L,\Lambda\}\), we have
\begin{equation}
\iota^{-1}(K_{\rmf}) = C_{M},\quad
\rho^{-1}(C_{M}) = j^{-1}(K_{\rmf}) = U^{1}.
\end{equation}
\end{lemma}

\begin{proof}
It is enough to prove this prime by prime.
Let \(\ell\) be a prime.
We first prove the first equality.
In the chosen coordinates, \(M_{\ell}=\Z_{\ell}^{3}\).
Let \(h\in \bfH(\Q_{\ell})\).
By \eqref{eq:PGLnZ}, \(\iota(h)\in \bfPGL_{3}(\Z_{\ell})\) if and only if \(ch\in \GL_{3}(\Z_{\ell})\) for some \(c\in \Q_{\ell}^{\times}\).
Since \(\det h=1\), \(c^{3}\in \Z_{\ell}^{\times}\), so \(c\in \Z_{\ell}^{\times}\), hence \(h\in \GL_{3}(\Z_{\ell})\).
Conversely, if \(h\in \GL_{3}(\Z_{\ell})\), \(\bbP(h)\in \bfPGL_{3}(\Z_{\ell})\).
Thus, \(\iota(h)\in \bfPGL_{3}(\Z_{\ell})\) if and only if \(h\in \GL_{3}(\Z_{\ell})\), if and only if \(hM_{\ell}=M_{\ell}\), i.e., \(\iota^{-1}(\bfPGL_{3}(\Z_{\ell}))=\rmop{Stab}_{\bfH(\Q_{\ell})}(M_{\ell})\).

We now prove the second equality.
We must prove that for every prime \(\ell\),
\[
\{s\in \bfS(\Q_{\ell})\mid \rho(s)M_{\ell}=M_{\ell}\}=\calO_{\ell}^{1}.
\]
(\(\supset\)).
If \(s\in \calO_{\ell}^{1}\), conjugation by \(s\) preserves \(\calO_{\ell}\), since \(s\) is a unit of \(\calO_{\ell}\), so conjugation by \(s\) preserves \(L_{\ell}=\ker(\trd|_{\calO_{\ell}})\) and \(\Lambda_{\ell}=T(\calO_{\ell})\) since reduced trace is invariant under conjugation and the map \(T\) commutes with conjugation.
Thus, \(\rho(\calO_{\ell}^{1})\) stabilizes the lattice \(M_{\ell}\) for either choice of \(M\in\{L,\Lambda\}\).

(\(\subset\)).
Now let \(s\in \bfS(\Q_{\ell})\) be such that \(\rho(s)\) stabilizes \(M_{\ell}\).
If \(\ell=p\), it follows from \cite[\S 6.4]{MaclachlanReid2003} that \(\bfS(\Q_{p})=\calO_{p}^{1}\) so we are done.
In the split case, if \(\ell\neq p\), we have an isomorphism \(\bfS(\Q_{\ell})\cong \SL_{2}(\Q_{\ell})\) and \(\calO_{\ell}\) is isomorphic to \(\rmop{M}_{2}(\Z_{\ell})\).
Using elementary divisors over \(\Z_{\ell}\), we can write \(s=u\rmop{diag}(\ell^{a},\ell^{-a})v\) with \(u,v\in \GL_{2}(\Z_{\ell})\) and \(a\in \Z\).
Conjugation by \(u\) and \(v\) preserves \(\calO_{\ell}=\rmop{M}_{2}(\Z_{\ell})\), so as before, conjugation by \(u\) and by \(v\) must preserve \(M_{\ell}\).
It follows that \(\rmop{diag}(\ell^{a},\ell^{-a})\) must also preserve \(M_{\ell}\).
The local description of \(L_{\ell}\) and \(\Lambda_{\ell}\) is immediate:
\[
L_{\ell}=\left\{\m{x & y \\ z & -x} \mid x,y,z\in \Z_{\ell}\right\},\quad
\Lambda_{\ell}=\left\{\m{x & 2y \\ 2z & -x} \mid x,y,z\in \Z_{\ell}\right\}.
\]
Then, it is easy to see the conjugation action as scaling the off-diagonal entries \((1,2)\) by \(\ell^{2a}\) and \((2,1)\) by \(\ell^{-2a}\).
Stabilizing the lattice forces \(a=0\) in either case, whence \(s\in \calO_{\ell}^{1}\).
It follows that, globally,  \(\rho^{-1}(C_{M})=\prod_{\ell} \calO_{\ell}^{1} = U^{1}\), and \(j^{-1}(K_{\rmf})=\rho^{-1}(\iota^{-1}(K_{\rmf}))=\rho^{-1}(C_{M})\).
\end{proof}

We thus have a diagram
\[
\begin{tikzcd}
\bfS(\bbA) \arrow[bend left]{rr}{j=\iota\circ \rho}\rar{\rho} & \bfH(\bbA) \rar{\iota} & \bfG(\bbA)\\
\bfS(\R)\times U^{1} \rar[swap]{\rho}\uar[hookrightarrow] & \bfH(\R) \times C_{M} \rar[swap]{\iota} \uar[hookrightarrow]& \bfG(\R) \times K_{\rmf}\uar[hookrightarrow]
\end{tikzcd}
\]
and in the finite adeles, we have a diagram of compact open subgroups.
The integral levels are as follows:
\[
\begin{tikzcd}
U^{1} =\calS(\hat{\Z})\dar[hookrightarrow]\rar{\rho} & C_{M}\dar[hookrightarrow]\rar{\iota} & K_{\rmf}=\bfPGL_{3}(\hat{\Z})\dar[hookrightarrow]\\
\bfS(\bbA_{\rmf}) \rar[swap]{\rho} & \bfH(\bbA_{\rmf}) \rar[swap]{\iota} & \bfG(\bbA_{\rmf})
\end{tikzcd}
\]
and \(W=\bfS(\R)\times U^{1}=j^{-1}(\bfG(\R)\times K_{\rmf})\) is a compact open subgroup of \(\bfS(\bbA)\).

\subsection{The homogeneous set}

We follow \cite[\S\S 1.2, 1.3]{EinsiedlerMargulisMohammadiVenkatesh2020}.
We now fix a relatively compact open set \(\Omega_{\infty}\subset \bfG(\R)\) containing \(\bbP(K_{\infty})\), and define the open subset \(\Omega_{0}=\Omega_{\infty} \times K_{\rmf} \subset \bfG(\bbA)\) containing the identity and has compact closure.
Volumes shall be considered relative to \(\Omega_{0}\) but changing \(\Omega_{0}\) keeps quantities comparable, as explained in \cite{EinsiedlerMargulisMohammadiVenkatesh2020}.

\begin{definition}
\label{definition:period-datum}
A \emph{datum} is a triple \(\scrD=(\bfT,j_{\bfT},a)\), where
\begin{enumerate}[auto]
\item \(\bfT\) is a connected semisimple algebraic \(\Q\)-group such that \(\bfT(\Q)\bs \bfT(\bbA)\) has finite volume;
\item The map \(j_{\bfT}:\bfT \to \bfG\) is an algebraic \(\Q\)-homomorphism with central finite kernel;
\item \(a\in \bfG(\bbA)\) is an element.
\end{enumerate}
We put \([\bfT]=\bfT(\Q)\bs\bfT(\bbA)\).
The datum defines the map
\[
\Phi_{\scrD}:[\bfT] \to X_{\bbA},\quad
\bfT(\Q)t \mapsto \bfG(\Q) j_{\bfT}(t) a,
\]
along with the \emph{homogeneous set} and its \emph{algebraic homogeneous measure}
\[
Y_{\scrD} = \Phi_{\scrD}([\bfT]),\quad
\mu_{\scrD}=(\Phi_{\scrD})_{*} m_{\bfT}.
\]
We say that \(\scrD\) or \(Y_{\scrD}\) is \emph{maximal} if \(j_{\bfT}(\bfT)_{\bar{\Q}}\subset \bfG_{\bar{\Q}}\) is a maximal proper algebraic subgroup.
We have an associated group
\[
H_{\scrD}=a^{-1} j_{\bfT}(\bfT(\bbA))a\subset \bfG(\bbA)
\]
which acts on \(Y_{\scrD}\) by right multiplication, and we let \(m_{\scrD}\) be the Haar measure on \(H_{\scrD}\) which projects to \(\mu_{\scrD}\) under the orbit map, i.e., the quotient measure on \(\rmop{Stab}_{H_{\scrD}}(\bfG(\Q)a)\bs H_{\scrD}\cong Y_{\scrD}\) is \(\mu_{\scrD}\).

The \emph{volume} of \(Y_{\scrD}\) relative to \(\Omega_{0}\) is defined as
\begin{equation}
\label{eq:hom-set-volume}
\rmop{vol}(Y_{\scrD}) = m_{\scrD}(H_{\scrD}\cap \Omega_{0})^{-1}.
\end{equation}
\end{definition}

For \(t\in \bfT(\bbA)\), we can write \(j_{\bfT}(t)a=a(a^{-1}j_{\bfT}(t)a)\).
It follows that
\begin{align*}
Y_{\scrD}
& = \Phi_{\scrD}(\bfT(\Q)\bs\bfT(\bbA))\\
& = \{\bfG(\Q)j_{\bfT}(t)a\mid \bfT(\Q)t\in \bfT(\Q)\bs\bfT(\bbA)\}\\
& = \{\bfG(\Q)ah\mid h\in H_{\scrD}\}\\
& = \bfG(\Q)aH_{\scrD}.
\end{align*}
If we put \(x_{0}=\bfG(\Q)a\in X_{\bbA}\), then \(Y_{\scrD}=x_{0}H_{\scrD}\), so \(Y_{\scrD}\) is the \(H_{\scrD}\)-orbit of \(x_{0}\) for right multiplication.
The group \(\bfG(\bbA)\) acts on \(X_{\bbA}\) by right translation: \((\bfG(\Q)g)\cdot h=\bfG(\Q)gh\).
Restricted to \(H_{\scrD}\), this action preserves \(Y_{\scrD}\) since the latter is an orbit.
Define the homomorphism
\[
\alpha_{\scrD}:\bfT(\bbA) \to H_{\scrD},\quad
s\mapsto a^{-1} j_{\bfT}(s)a.
\]
Then, \(H_{\scrD}=\alpha_{\scrD}(\bfT(\bbA))\), and \(\alpha_{\scrD}\) makes \(\Phi_{\scrD}\) equivariant for the right action,  because for \(s\in \bfT(\bbA)\), we have \(j_{\bfT}(s)a = a \alpha_{\scrD}(s)\).
Hence, for \(s,t\in \bfT(\bbA)\),
\[
\Phi_{\scrD}(\bfT(\Q)ts)=\bfG(\Q)j_{\bfT}(ts)a
=\Phi_{\scrD}(\bfT(\Q)t)\alpha_{\scrD}(s).
\]

In order to use the results of \cite{EinsiedlerMargulisMohammadiVenkatesh2020}, we need maximality of the datum, which we now prove.
To prove maximality, we use results on the classification of maximal subgroups of classical groups (see \cite[Chapter 18]{MalleTesterman2011}).
Maximal subgroups are obtained as normalizers of specific groups;
from the classification of maximal subgroups of \(\bfSL_{3}\), we can then show that the corresponding subgroup in \(\bfG\) is maximal.
The maximal subgroup classification of \cite[Chapter 18]{MalleTesterman2011} is formulated for Zariski-closed subgroups of the groups of \(k\)-points for an algebraically closed field \(k\).
In characteristic zero, this is not an obstacle to our group scheme theoretic framework, see the discussion in \cite[\S\S 3.f and 3.g]{Milne2017}.
In particular, \cite[Proposition 3.16, Theorem 3.23, and Corollaries 3.26, 3.27]{Milne2017} suffice for our purposes, to check for maximality of proper algebraic subgroups (in particular, \cite[Corollary 3.27]{Milne2017} concerns normalizers).

\begin{lemma}
\label{lemma:geometric-maximality}
The subgroup \(\iota(\bfH)_{\bar{\Q}}\) is a maximal proper algebraic subgroup of \(\bfG_{\bar{\Q}}\).
\end{lemma}

\begin{proof}
Crucially, there is no assumption on connectedness, so we must show that \(\iota(\bfH)_{\bar{\Q}}\) is maximal among all proper algebraic subgroups of \(\bfG_{\bar{\Q}}\), connected or not.

In this proof, let \(k=\bar{\Q}\), \(V_{k}=V\otimes_{\Q}k\), \(\bfH_{0} = \bfSO(V_{k},q)=\bfH_{k}\), and let \(\pi:\bfSL(V_{k}) \to \bfPGL(V_{k})\) be the quotient map (the composite \(\bfSL_{3} \into \bfGL_{3} \to \bfPGL_{3}\)).
This is a central isogeny with kernel \(\bfmu_{3}\).
The coordinates chosen identify \(\bfSL(V_{k})\) with \(\bfSL_{3,k}\), \(\bfPGL(V_{k})\) with \(\bfG_{k}\), and \(\pi(\bfH_{0})\) with \(\iota(\bfH)_{k}\).
Since algebraic subgroups are determined by their \(k\)-points, we can apply results of \cite{MalleTesterman2011}.

First, we observe that \(V_{k}\) is a simple \(\bfH_{0}\)-module.
By \cite[Theorem 5.5.13]{GoodmanWallach2009}, the defining representation of \(\SO(\C^{3},B)\) on \(\C^{3}\) is irreducible for any nondegenerate symmetric bilinear form \(B\).
After a change of basis, this implies in particular that the same holds for \(\SO(V_{\C},q)\) on \(V_{\C}\).
Fix an embedding \(k\into \C\).
If \(W\subset V_{k}\) were a proper nontrivial \(\bfH_{0}\)-invariant subspace, then \(W\otimes_{k}\C\) would be a proper nontrivial subspace of \(V_{\C}\) invariant under \(\bfH_{0}(\C)=\SO(V_{\C},q)\), a contradiction.
Thus, \(V_{k}\) is a simple \(\bfH_{0}\)-module.

Next, we claim that
\[
N_{\bfSL(V_{k})}(\bfH_{0}) = \bfmu_{3} \bfH_{0} = \pi^{-1}(\pi(\bfH_{0})).
\]
This equality identifies \(\bfmu_{3}\bfH_{0}\) with the normalizer of \(\bfH_{0}\), a subgroup in the class \(\calC_{6}\) from \cite[Table 18.2]{MalleTesterman2011}.
The second equality holds because \(\ker\pi=\bfmu_{3}\).
For the first equality, we first have \(\bfmu_{3}\bfH_{0}\subset N_{\bfSL(V_{k})}(\bfH_{0})\) because \(\bfmu_{3}\) is central.
For the reverse inclusion, let \(g\in N_{\bfSL(V_{k})}(\bfH_{0})(k)\) and let \(Q\) be a Gram matrix of \(q\) on \(V_{k}\).
Since \(g\) normalizes \(\bfH_{0}\), the bilinear form \(b_{g}(v,w)=v^{T}(g^{T}Qg)w\) with Gram matrix \(g^{T}Qg\) is also \(\bfH_{0}\)-invariant, so that \(h^{T}(g^{T}Qg)h=g^{T}Qg\) for every \(h\in \bfH_{0}(k)\).
Let \(C=Q^{-1}(g^{T}Qg)\).
We claim that \(C\) commutes with every \(h\in \bfH_{0}(k)\).
Since \(h^{T}Qh=Q\), we have \(h^{-1}=Q^{-1}h^{T}Q\), so
\[
h^{-1}Ch=(Q^{-1}h^{T}Q)(Q^{-1}(g^{T}Qg))h
=Q^{-1}h^{T}g^{T}Qgh=Q^{-1}g^{T}Qg=C.
\]
Thus, \(C\in \End_{\bfH_{0}}(V_{k})\).
Since \(k\) is algebraically closed, by Schur's lemma, \(C=cI_{3}\) for some \(c\in k^{\times}\), that is, \(g^{T}Qg=cQ\) for some \(c\in k^{\times}\).
Taking determinants, we get \(c^{3}=1\).
Let \(\zeta=c^{2}\in \bfmu_{3}(k)\);
then \(\zeta^{2}=c\), so
\[
(\zeta^{-1}g)^{T}Q(\zeta^{-1}g)=\zeta^{-2}g^{T}Qg=\zeta^{-2}cQ=Q,\quad
\det(\zeta^{-1}g)=\zeta^{-3}\det g =1.
\]
Thus, \(\zeta^{-1}g\in \bfH_{0}(k)\).
We conclude that \(g=\zeta h_{0}\in \bfmu_{3}(k) \bfH_{0}(k)\), which establishes the equality.

We can now apply \cite[Proposition 18.13, Table 18.2, class \(\calC_{6}\)]{MalleTesterman2011} for the special case of \(\bfSL(V_{k})\cong \bfSL_{3,k}\) and the orthogonal group of class \(\calC_{6}\).
These results are due to Liebeck and Seitz \cite{LiebeckSeitz1998}.
This implies that the normalizer \(N_{\bfSL(V_{k})}(\bfH_{0})\) is a maximal proper closed subgroup of \(\bfSL(V_{k})\).
We can now conclude:
if \(\bfJ\) is a closed algebraic subgroup such that \(\pi(\bfH_{0})\subset \bfJ \subset \bfPGL(V_{k})\), then
\[
N_{\bfSL(V_{k})}(\bfH_{0}) = \pi^{-1}(\pi(\bfH_{0})) \subset \pi^{-1}(\bfJ) \subset \bfSL(V_{k}).
\]
Maximality then forces \(\pi^{-1}(\bfJ)\) to be either \(\bfSL(V_{k})\) or \(N_{\bfSL(V_{k})}(\bfH_{0})\).
Since \(\pi\) is surjective, \(\bfJ=\pi(\bfH_{0})\), or \(\bfJ=\bfPGL(V_{k})\), completing the proof.
\end{proof}

The datum attached to the lattice \(M\) is \(\scrD_{p}=\scrD_{p,M}=(\bfS,j_{M},\omega_{p,M})\), see \eqref{eq:future-datum}.
Correspondingly, we write
\begin{equation}
\label{eq:definition-special-homogeneous-set}
\left\{
\begin{alignedat}{2}
\Phi_{p}:\bfS(\Q)\bs\bfS(\bbA) &\to X_{\bbA},&\qquad
\bfS(\Q)s&\mapsto \bfG(\Q) j(s)\omega_{p},\\
Y_{p}=Y_{p,M} & = \Phi_{p}(\bfS(\Q)\bs\bfS(\bbA)),&\qquad
\mu_{p}=\mu_{p,M} & = (\Phi_{p})_{*}m_{\bfS},
\end{alignedat}
\right.
\end{equation}
and set \(Z=\ker (\bfS(\bbA) \xto{j}\bfG(\bbA)) =\bfmu_{2}(\bbA) = \prod_{v} \{\pm 1\}\), along with
\begin{equation}
\left\{
\begin{alignedat}{2}
\pi:\bfS(\bbA) & \to \bfS(\bbA)/Z,\\
\alpha_{p}:\bfS(\bbA) &\to H_{\scrD_{p}},&\quad
\alpha_{p}(s)&=\omega_{p}^{-1}j(s)\omega_{p}.
\end{alignedat}
\right.
\end{equation}
(Note that \(\alpha_{p}=\alpha_{\scrD_{p}}\) in the notation of Definition \ref{definition:period-datum}.)

To ease notation, we introduce the following notation attached to the group \(\bfS\).
We put
\[
G^{1} = \bfS(\bbA),\quad
\Gamma^{1} = \bfS(\Q),\quad
[\bfS]=\bfS(\Q)\bs \bfS(\bbA)=\Gamma^{1}\bs G^{1}.
\]
We emphasize that even though \(\bfG=\bfPGL_{3}\), \(G^{1}\) corresponds to the norm-one group \(\bfS\).
For quotients modulo centers, we also put
\[
\bar{G^{1}} = G^{1}/Z,\quad
\bar{\Gamma^{1}} = \Gamma^{1}Z/Z\cong \Gamma^{1}/\{\pm 1\}.
\]
Then,
\[
\pi(\Gamma^{1})= \Gamma^{1}Z/Z\cong \Gamma^{1}/\{\pm 1\}.
\]
The subgroup \(\Gamma^{1}\) is discrete in \(G^{1}\), and \(\Gamma^{1}\cap Z=\{\pm 1\}\).
The quotient \([\bfS]\) is compact by Godement's compactness criterion \cite[Theorem 5.6]{Borel1963}:
\(\bfS\) is semisimple, hence has no nontrivial \(\Q\)-characters, and \(\bfS(\Q)\) contains no nontrivial unipotent elements since \(B\) is a division algebra.
Let \(m_{\bfS}\) denote its invariant probability measure.
The kernel \(Z\) is compact, central, and contained in \(W\);
\(\alpha_{p}\) is surjective and has kernel \(\ker(\alpha_{p})=\ker(j)=Z\), so that \(\alpha_{p}=\bar{\alpha}_{p}\circ \pi\) for a bijection \(\bar{\alpha}_{p}:\bar{G^{1}}\to H_{\scrD_{p}}\), and \(\bar{\alpha}_{p}(\bar{\Gamma^{1}})=\alpha_{p}(\Gamma^{1})=\omega_{p}^{-1} j(\Gamma^{1})\omega_{p}\).
Since \(Z\) is central, \(\Gamma^{1} sZ=\Gamma^{1}Zs\) for \(s\in G^{1}\), and we can identify \(\bar{\Gamma^{1}}\bs \bar{G^{1}}=[\bfS]/Z\), the quotient of \([\bfS]\) by the right action of \(Z\).
Moreover,
\begin{equation}
\label{eq:compatibility-homogeneous}
\bfG(\Q) \omega_{p} \alpha_{p}(s) = \bfG(\Q) j(s) \omega_{p} = \Phi_{p}(\Gamma^{1}s)\quad (s\in G^{1}).
\end{equation}

We proceed to analyzing the homogeneous set \(Y_{p}\).

\begin{proposition}
\label{proposition:structure-homogeneous-set}
Let \(M\in \{L,\Lambda\}\).
\begin{enumerate}[auto]
\item \(\scrD_{p}\) is a maximal datum with \(\bfS\) simply connected;
\(Y_{p}\) is compact, and \(\mu_{p}\) is a probability measure.

\item \(j( G^{1} )\cap \bfG(\Q) = j( \Gamma^{1} )\).

\item The subgroup \(j( G^{1} )\) is closed in \(\bfG(\bbA)\), and \(\bar{\alpha}_{p}: \bar{G^{1}} \stosim H_{\scrD_{p}}\) is an isomorphism of topological groups.
In particular, the subgroup \(H_{\scrD_{p}}\) is also closed in \(\bfG(\bbA)\), and \(\alpha_{p}: G^{1}  \to H_{\scrD_{p}}\) is continuous, open, and proper.

\item
The stabilizer of \(\bfG(\Q)\omega_{p}\) in \(H_{\scrD_{p}}\) is the discrete subgroup \(\alpha_{p}( \Gamma^{1} )=\omega_{p}^{-1}j( \Gamma^{1} )\omega_{p}\), and the diagram
\[
\begin{tikzcd}
G^{1} \dar \rar{\pi} &\bar{G^{1}} \dar \rar{\bar{\alpha}_{p}}[swap]{\simeq} & H_{\scrD_{p}}\dar\arrow{rd}{h\mapsto \bfG(\Q)\omega_{p}h}\\
{}[\bfS] =\Gamma^{1}\bs G^{1} \rar[twoheadrightarrow] & \bar{\Gamma^{1}}\bs \bar{G^{1}} \rar{\simeq} & \alpha_{p}(\Gamma^{1})\bs H_{\scrD_{p}} \rar{\simeq} & Y_{p}
\end{tikzcd}
\]
commutes, where \(\simeq\) here denotes homeomorphisms.
The vertical arrows are quotient maps, and the three bottom arrows are
\[
\Gamma^{1}s \mapsto \bar{\Gamma^{1}}\pi(s),\quad
\bar{\Gamma^{1}}x \mapsto \alpha_{p}(\Gamma^{1})\bar{\alpha}_{p}(x),\quad
\alpha_{p}(\Gamma^{1})h\mapsto \bfG(\Q) \omega_{p}h,
\]
respectively.
The map \([\bfS] \onto \bar{\Gamma^{1}}\bs\bar{G^{1}}\) is surjective, with fibers the \(Z\)-orbits \(\Gamma^{1}sZ\).
The bottom composition is \(\Phi_{p}:[\bfS] \to Y_{p}\subset X_{\bbA}\).
In particular, \(\Phi_{p}\) induces a homeomorphism
\[
\bar{\Phi}_{p}:\bar{\Gamma^{1}}\bs \bar{G^{1}} \to Y_{p},\quad
\bar{\Phi}_{p}\bigl(\bar{\Gamma^{1}}\pi(s)\bigr)=\bfG(\Q) j(s) \omega_{p}.
\]
\end{enumerate}
\end{proposition}

\begin{proof}
(1)
The group \(\bfS\) is connected, semisimple, and simply connected, and \(\ker j=\bfmu_{2}\) is finite and central.
Since \(\rho\) is surjective, \(j(\bfS)=\iota(\bfH)\), which is maximal by Lemma \ref{lemma:geometric-maximality}.
Since \([\bfS]\) is compact and \(\Phi_{p}\) continuous, it follows that \(Y_{p}\) is compact.

(2)
We first establish
\begin{equation}
\tag{i}\label{eq:4.6.2.i}
\iota(\bfH(\bbA))\cap \bfG(\Q)=\iota(\bfH(\Q)).
\end{equation}
Let \(x\in \iota(\bfH(\bbA))\cap \bfG(\Q)\).
Since \(\iota\) is a closed immersion over \(\Q\), the image \(\iota(\bfH)\) is defined by an ideal \(I\subset \Q[\bfG]\) (\cite[Chapter 2, \S 3, Proposition 3.20]{Liu2002}).
Every \(f\in I\) satisfies \(f(x)=0\) in \(\bbA\), hence already in \(\Q\), since \(\Q\into\bbA\) is injective, thus \(x\in \iota(\bfH)(\Q)=\iota(\bfH(\Q))\) since \(\iota\) is a closed immersion.
The reverse inclusion \(\iota(\bfH(\bbA))\cap \bfG(\Q)\supset\iota(\bfH(\Q))\) holds because \(\iota\) is defined over \(\Q\).

We now establish
\begin{equation}
\tag{ii}\label{eq:4.6.2.ii}
\rho( G^{1} )\cap \bfH(\Q)=\rho( \Gamma^{1} ).
\end{equation}
Let \(h\in \rho( G^{1} )\cap \bfH(\Q)\).
since \(\rho=\Psi|_{\bfS}\), we can write \(h=\Psi(b)\) with \(b\in B^{\times}\).
We use the criterion that \(\Psi(b)\in \rho(\bfS(F))\) if and only if \(\nrd(b)\in F^{\times 2}\) for any field extension \(F/\Q\).
At every place \(v\), choose \(s_{v}\in \bfS(\Q_{v})\) such that \(\rho(s_{v})=h_{v}=\Psi(s_{v})\).
Since \(\ker\Psi(\Q_{v})=\Q_{v}^{\times}\), we have \(b=c_{v}s_{v}\) for some \(c_{v}\in \Q_{v}^{\times}\).
Since \(s_{v}\in \bfS(\Q_{v})\), we have \(\nrd(b)=c_{v}^{2}\) at every place.
Since \(\nrd(b)\) is a positive rational number which is a local square at every place, it must be a global square (the valuation at every prime is even).
Hence, there is \(c\in \Q^{\times}\) such that \(\nrd(b)=c^{2}\).
Then, \(c^{-1}b\in  \Gamma^{1} \), and \(\rho(c^{-1}b)=\Psi(c^{-1}b)=\Psi(b)=h\), since \(\Psi\) is conjugation and \(c\) is a scalar.
Thus, we have proved \(\rho( G^{1} )\cap \bfH(\Q)\subset \rho( \Gamma^{1} )\).
The reverse inclusion holds because \(\rho\) is defined over \(\Q\).

Finally, we can prove \(j( G^{1} )\cap \bfG(\Q)=j( \Gamma^{1} )\).
Since \(\iota(\rho( G^{1} ))\subset \iota(\bfH(\bbA))\), intersecting with \(\bfG(\Q)\), we have, since \(\iota\) is injective on \(\bfH(\bbA)\),
\begin{alignat*}{2}
 \iota(\rho( G^{1} )) \cap \bfG(\Q)
& = \iota(\rho( G^{1} ))\cap \iota(\bfH(\bbA))\cap \bfG(\Q)\\
& = \iota(\rho( G^{1} ))\cap \iota(\bfH(\Q)) &&\qquad (\text{by \eqref{eq:4.6.2.i}})\\
& = \iota\bigl(\rho( G^{1} )\cap \bfH(\Q)\bigr)\\
& = \iota(\rho( \Gamma^{1} )) &&\qquad (\text{by \eqref{eq:4.6.2.ii}}).
\end{alignat*}
Thus, \(j( G^{1} )\cap \bfG(\Q) = j( \Gamma^{1} )\).

(3)
Let \(\Omega=\bfG(\R)\times K_{\rmf}\).
This is an open subgroup of \(\bfG(\bbA)\).
By Lemma \ref{lemma:integral-structures}, \(j^{-1}(\Omega)=W=\bfS(\R)\times U^{1}\).
Since \(W\) is compact, \(j(W)=j( G^{1} )\cap \Omega\) is compact, hence closed.
Let \(x\in j( G^{1} )\cap g\Omega\).
Since \(j( G^{1} )\) is a subgroup, \(j( G^{1} )\cap g\Omega=x(j( G^{1} )\cap\Omega)=xj(W)\), which is compact, hence closed in \(g\Omega\).
Since the cosets \(g\Omega\) form a partition of \(\bfG(\bbA)\) made of open sets, it follows that \(j( G^{1} )\) is closed.

The map \(W/Z\to j(W)\) is a continuous bijection from a compact space to a Hausdorff space, so it is a homeomorphism.
The map \(\bar{j}: \bar{G^{1}} \to j( G^{1} )\) is a continuous bijection.
Restricted to \(W/Z\), it is a homeomorphism, and restricted to the coset \(\pi(s)(W/Z)\), it is also a homeomorphism.
The cosets \(\pi(s)(W/Z)\), \(s\in G^{1}\), form open covers of \( \bar{G^{1}} \), and their images  \(j(s)j(W)\) form an open cover of \(j( G^{1} )\).
Thus, \(\bar{j}: \bar{G^{1}}\stosim j( G^{1} )\) is an isomorphism of topological groups.
Since \(\alpha_{p}\) is obtained by conjugation by \(\omega_{p}\), we conclude that \(\bar{\alpha}_{p}: \bar{G^{1}} \to H_{\scrD_{p}}\) is an isomorphism of topological groups.
Since \(\pi\) is open and \(Z\) is compact, \(\alpha_{p}=\bar{\alpha}_{p}\circ\pi\) is open and proper.

(4)
Recall that \(\Phi_{p}( \Gamma^{1} s)=\bfG(\Q)j(s)\omega_{p}=\bfG(\Q)\omega_{p}\alpha_{p}(s)\) for \(s\in  G^{1} \) by \eqref{eq:compatibility-homogeneous}.
For \(s\in G^{1}\), by part (2),
\[
\bfG(\Q)\omega_{p}\alpha_{p}(s)=\bfG(\Q)\omega_{p}
\;\Leftrightarrow\;
j(s)\in \bfG(\Q) \cap j(G^{1})=j(\Gamma^{1})
\;\Leftrightarrow\;
s\in \Gamma^{1}Z,
\]
since \(\ker j=Z\).
Hence, the stabilizer of \(\bfG(\Q)\omega_{p}\) is precisely \(\alpha_{p}(\Gamma^{1}Z)=\alpha_{p}(\Gamma^{1}\subset \omega_{p}^{-1}\bfG(\Q)\omega_{p}\), so it is a discrete subgroup since \(\bfG(\Q)\) is discrete in \(\bfG(\bbA)\).

Since \(Z\) is central, \(\bar{\Gamma^{1}}\pi(t)=\bar{\Gamma^{1}}\pi(s)\) if and only if \(t\in \Gamma^{1}sZ\), so the fibers of the map \([\bfS]\to \bar{\Gamma^{1}}\bs\bar{G^{1}}\) are the \(Z\)-orbits \(\Gamma^{1}sZ\).
The arrow \(\bar{\Gamma^{1}}\bs\bar{G^{1}}\to \alpha_{p}(\Gamma^{1})\bs H_{\scrD_{p}}\) is a homeomorphism:
it is induced by \(\bar{\alpha}_{p}\), an isomorphism of topological groups carrying \(\bar{\Gamma^{1}}\) onto \(\alpha_{p}(\Gamma^{1})\).
The final arrow \(\alpha_{p}(\Gamma^{1})\bs H_{\scrD_{p}}\to Y_{p}\) is a well-defined homeomorphism, since it is induced by the orbit map, the stabilizer is \(\alpha_{p}(\Gamma^{1})\), and the source is compact, as the image of \([\bfS]\) with Hausdorff target.
By definition, \(Y_{p}=\Phi_{p}([\bfS])=\bfG(\Q)\omega_{p}\alpha_{p}(G^{1})=\bfG(\Q)\omega_{p}H_{\scrD_{p}}\).

Commutativity of the diagram is easily checked by inspection of the explicit maps in the bottom row.
\end{proof}

We now need to compute the volume \(\rmop{vol}(Y_{p})\) and the pushforward \(\rmop{pr}_{*}\mu_{p}\).

\subsection{The volume}

Let \(ds\) be the Haar measure on \(G^{1}\) whose quotient measure on \([\bfS]\) is \(m_{\bfS}\).
We need to compute the volume of \(Y_{p}\) using \(W=\bfS(\R)\times U^{1}\).

\begin{lemma}
\label{lemma:fiber-measure}
For \(s=(s_{\infty},s_{\rmf})\in  G^{1} \), the map
\[
\rmop{cl}:[\bfS] \to \scrC_{\bfS},\quad
\Gamma^{1}s \mapsto \Gamma^{1}s_{\rmf}U^{1}=[s_{\rmf}],
\]
induces a bijection \([\bfS]/W\to \scrC_{\bfS}\).
The fiber above the double coset \(\frakc=[t_{\rmf}]=\Gamma^{1}t_{\rmf}U^{1}\in \scrC_{\bfS}\) is the compact open \(W\)-orbit
\[
\scrT_{\frakc} = \cl^{-1}(\frakc)=\Gamma^{1}\bs \Gamma^{1}(1,t_{\rmf})W\subset[\bfS],
\]
and
\begin{equation}
\label{eq:orbit-volume}
ds(W)=\frac{1}{M_{p}},\quad
m_{\bfS}(\scrT_{\frakc})=\frac{1}{M_{p}|\calO_{t_{\rmf}}^{\times}|} .
\end{equation}
\end{lemma}

\begin{proof}
The map \(\rmop{cl}\) is well-defined and constant on \(W\)-orbits:
for \(g\in  G^{1} \) and \(w=(w_{\infty},u)\in W\), \(\rmop{cl}(\Gamma^{1}gw)=\Gamma^{1} g_{\rmf}uU^{1}=\Gamma^{1}g_{\rmf}U^{1}=\rmop{cl}(\Gamma^{1}g)\).
Conversely, if \(\rmop{cl}(\Gamma^{1}g^{\prime})=\rmop{cl}(\Gamma^{1}g)\) with \(g,g^{\prime}\in  G^{1} \), then \(g_{\rmf}^{\prime}=\gamma_{\rmf}g_{\rmf}u\) for some \(\gamma\in \Gamma^{1}\) and \(u\in U^{1}\).
Then, since \(\Gamma^{1}\) is embedded diagonally,
\[
g^{\prime}=(g_{\infty}^{\prime},\gamma g_{\rmf}u)=\gamma g ((\gamma g_{\infty})^{-1}g_{\infty}^{\prime},u),
\]
and \(w=((\gamma g_{\infty})^{-1}g_{\infty}^{\prime},u)\in \bfS(\R)\times U^{1}=W\), so \(g^{\prime}=\gamma g w\).
Thus, \(\Gamma^{1}g^{\prime}\) and \(\Gamma^{1}g\) lie in the same right \(W\)-orbit, hence the fibers of \(\rmop{cl}\) are exactly the \(W\)-orbits, and \(\cl\) induces the asserted bijection, surjective because \(\rmop{cl}(\Gamma^{1}(1,s_{\rmf}))=[s_{\rmf}]\).
For \(\frakc=\Gamma^{1}t_{\rmf}U^{1}\), we clearly have \(\scrT_{\frakc}=\Gamma^{1}\bs \Gamma^{1}(1,t_{\rmf})W\), which is compact open, since \(W\) is compact open and the projection \( G^{1} \to [\bfS]\) is open.

Fix \(\frakc=[t_{\rmf}]\in \scrC_{\bfS}\), and put \(g=(1,t_{\rmf})\), and \(\scrT=\scrT_{\frakc}\).
By Theorem \ref{theorem:yu-norm-one-correspondence} and \eqref{eq:isomorphism-conjugation-endomorphism-ring}, \(\Gamma^{1}\cap t_{\rmf}U^{1}t_{\rmf}^{-1}=\calO_{t_{\rmf}}^{\times}\cong\Aut(E)\), with \([E]=\Theta^{-1}([t_{\rmf}])\), which is finite.
Since \(gWg^{-1}=\bfS(\R)\times t_{\rmf}U^{1}t_{\rmf}^{-1}\) and \(\Gamma^{1}\) is embedded diagonally, the real component imposes no condition, and \eqref{eq:isomorphism-conjugation-endomorphism-ring} gives
\[
\Gamma^{1} \cap gWg^{-1} = \Gamma^{1}\cap t_{\rmf}U^{1}t_{\rmf}^{-1}=\calO_{t_{\rmf}}^{\times},
\]
a finite set.
Let \(\psi=\mathbf{1}_{gW}\in C_{c}( G^{1} )\), since \(gW\) is compact open.
By \eqref{eq:Weil},
\[
ds(W)=ds(gW) = \int_{\Gamma^{1}\bs  G^{1} } \sum_{\gamma\in \Gamma^{1}} \psi(\gamma x) dm_{\bfS}(x)
=\int_{[\bfS]} F\, dm_{\bfS},
\]
where \(F(\Gamma^{1}x)=\#\{\gamma\in \Gamma^{1} \mid \gamma x\in gW\}\);
\(F\) is well defined on \([\bfS]\), since replacing \(x\) by \(\gamma_{0}x\) reindexes the count.
If \(\Gamma^{1}x\notin \scrT\), then \(x\notin \Gamma^{1}gW\), then \(\gamma x\notin gW\) for every \(\gamma\in \Gamma^{1}\), hence \(\psi(\gamma x)=0\) for every \(\gamma\in \Gamma^{1}\) in  this case, so \(F(\Gamma^{1}x)=0\) for \(\Gamma^{1}x\notin \scrT\).
Otherwise, \(\Gamma^{1}x\in \scrT\), so \(x\in \Gamma^{1}gW\), and we can write \(x=\gamma_{0}gw\) with \(\gamma_{0}\in \Gamma^{1}\) and \(w\in W\).
Then,
\[
\gamma x \in gW\Leftrightarrow
\gamma \gamma_{0} gw \in gW\Leftrightarrow
g^{-1}\gamma \gamma_{0}g\in W\Leftrightarrow
\gamma\gamma_{0}\in gWg^{-1}.
\]
Therefore, \(F(\Gamma^{1}x)=|\Gamma^{1}\cap gWg^{-1}|=|\calO_{t_{\rmf}}^{\times}|\) in this case.
We conclude that \(F=|\calO_{t_{\rmf}}^{\times}|\, \mathbf{1}_{\scrT}\), so \eqref{eq:Weil} gives
\[
ds(W) =|\calO_{t_{\rmf}}^{\times}|\cdot m_{\bfS}(\scrT).
\]

The fibres \(\scrT_{\frakc}\) form a finite partition of \([\bfS]\), and \(m_{\bfS}([\bfS])=1\).
By Theorem \ref{theorem:yu-norm-one-correspondence} and \eqref{eq:isomorphism-conjugation-endomorphism-ring},
\(|\calO_{t_{\rmf}}^{\times}|=|\Aut(E)|\) for the curve \(E\) corresponding to \([t_{\rmf}]\), so by \eqref{eq:mass-formula-Aut},
\[
1=\sum_{\frakc\in \scrC_{\bfS}} m_{\bfS}(\scrT_{\frakc})
=ds(W) \sum_{[t_{\rmf}]\in \scrC_{\bfS}}\frac{1}{|\calO_{t_{\rmf}}^{\times}|}
=ds(W) \cdot M_{p}.
\]
Thus,
\[
ds(W) = \frac{1}{M_{p}},\quad
m_{\bfS}(\scrT_{\frakc}) = \frac{ds(W)}{|\calO_{t_{\rmf}}^{\times}|}=\frac{1}{M_{p}\, |\calO_{t_{\rmf}}^{\times}|}.
\]
This completes the proof.
\end{proof}

We can now compute the volume of the homogeneous set.

\begin{proposition}
\label{proposition:homogeneous-set-volume}
The Haar measure \(m_{\scrD_{p}}\) on \(H_{\scrD_{p}}\) satisfies \(m_{\scrD_{p}}=\frac{1}{2}(\alpha_{p})_{*}ds\), and \(\alpha_{p}^{-1}(\Omega_{0})=W\).
As a consequence, for \(M\in\{L,\Lambda\}\), we have
\begin{equation}
\label{eq:volume-homogeneous-set}
v_{p} = \rmop{vol}(Y_{p,M}) = 2M_{p} = \frac{p-1}{12}.
\end{equation}
\end{proposition}

\begin{proof}
Let \(\nu=(\alpha_{p})_{*}ds\).
Then, \(\nu\) is a Haar measure on \(H_{\scrD_{p}}\).
By Proposition \ref{proposition:structure-homogeneous-set}, \(\alpha_{p}( \Gamma^{1} )\) is the stabilizer of \(\bfG(\Q)\omega_{p}\), and the orbit map induces a homeomorphism \(\alpha_{p}( \Gamma^{1} )\bs H_{\scrD_{p}} \to Y_{p}\).

Let \(f\in C_{c}(H_{\scrD_{p}})\).
By \eqref{eq:Weil},
\begin{align*}
\int_{H_{\scrD_{p}}} f \, d\nu
&= \int_{ G^{1} } f(\alpha_{p}(s))\, ds
= \int_{[\bfS]} \sum_{\gamma\in  \Gamma^{1} } f(\alpha_{p}(\gamma)\alpha_{p}(s))\, dm_{\bfS}\\
& = 2 \int_{[\bfS]} \sum_{\lambda \in \alpha_{p}( \Gamma^{1} )} f(\lambda \alpha_{p}(s)) \, dm_{\bfS}.
\end{align*}
The functions are functions of \( \Gamma^{1} s\in [\bfS]\), and the last equality holds because \(\alpha_{p}|_{ \Gamma^{1} }\) is surjective onto its image \(\alpha_{p}( \Gamma^{1} )\) with kernel \( \Gamma^{1} \cap Z = \{\pm 1\}\), so that each \(\lambda\) above is \(\alpha_{p}(\pm\gamma)\) for \(\gamma\in  \Gamma^{1} \).
Since \(\alpha_{p}( \Gamma^{1} )\) is the stabilizer of \(\bfG(\Q)\omega_{p}\), the function \(F(\bfG(\Q)\omega_{p}h)=\sum_{\lambda\in \alpha_{p}( \Gamma^{1} )}f(\lambda h)\) is well defined on \(Y_{p}\).
By \eqref{eq:compatibility-homogeneous},
\[
\sum_{\lambda}f(\lambda \alpha_{p}(s)) = F(\bfG(\Q)\omega_{p}\alpha_{p}(s)) = F(\Phi_{p}( \Gamma^{1} s)).
\]
We conclude that
\[
\int_{H_{\scrD_{p}}} f \, d\nu
=2\int_{[\bfS]} F(\Phi_{p}( \Gamma^{1} s)) dm_{\bfS}=2\int_{Y_{p}} F\, d\mu_{p},
\]
since \(\mu_{p}=(\Phi_{p})_{*}m_{\bfS}\).
Thus, the quotient measure of \(\nu\) by \(\alpha_{p}( \Gamma^{1} )\), in the sense of the Weil integration formula \eqref{eq:Weil}, transported to \(Y_{p}\) by the orbit map, is \(2\mu_{p}\).
So \(\frac{1}{2}\nu\) projects to \(\mu_{p}\), and \(m_{\scrD_{p}}\) is defined as the measure projecting to \(\mu_{p}\), we conclude that \(m_{\scrD_{p}}=\frac{1}{2}\nu=\frac{1}{2}(\alpha_{p})_{*}ds\), as desired.

We have \(\omega_{p}=(\bbP(\xi ),1_{\rmf})\).
For \(s=(s_{\infty},s_{\rmf})\),
\[
\alpha_{p}(s)=\omega_{p}^{-1}j(s)\omega_{p}=(\bbP(\xi ^{-1}\rho(s_{\infty})\xi ),j(s_{\rmf})).
\]
By \eqref{eq:kernels-compact-subgroup}, the archimedean component lies in \(\bbP(K_{\infty})\subset \Omega_{\infty}\), and the finite component lies in \(K_{\rmf}\) exactly when \(s_{\rmf}\in U^{1}\), by Lemma \ref{lemma:integral-structures}.
Since \(\rmop{im}(\alpha_{p})\subset H_{\scrD_{p}}\subset \bfG(\bbA)\) and \(\Omega_{0}=\Omega_{\infty}\times K_{\rmf}\subset \bfG(\bbA)\),
\[
\alpha_{p}^{-1}(H_{\scrD_{p}}\cap\Omega_{0})
=\{s\in  G^{1}  \mid \alpha_{p}(s) \in \Omega_{0}\}
=\bfS(\R)\times U^{1}=W.
\]
By \eqref{eq:hom-set-volume}, the definition of \((\alpha_{p})_{*}ds\), and \eqref{eq:orbit-volume},
\begin{align*}
\rmop{vol}(Y_{p,M})
& = \bigl(\tfrac{1}{2}(\alpha_{p})_{*} ds(H_{\scrD_{p}}\cap \Omega_{0})\bigr)^{-1}
 = \bigl(\tfrac{1}{2}ds\bigl(\alpha_{p}^{-1}(H_{\scrD_{p}}\cap \Omega_{0})\bigr)\bigr)^{-1}\\
 &= \bigl(\tfrac{1}{2}( ds(W)\bigr)^{-1}
= 2ds(W)^{-1}
 = 2M_{p}
 = \frac{p-1}{12}.
\end{align*}
This completes the proof.
\end{proof}

\subsection{Projecting to shapes}

\begin{proposition}
\label{proposition:final-projection}
Let \(M\in \{L,\Lambda\}\).
For every \(h=(h_{\infty},h_{\rmf})\in \bfH(\bbA)\),
\begin{equation}
\label{eq:projection-to-shapes}
\rmop{pr}\bigl(\bfG(\Q)\iota(h)\omega_{p}\bigr) = \rmop{sh}(h_{\rmf}\cdot M,q).
\end{equation}
The diagram \eqref{eq:goal-diagram} commutes, and
\begin{equation}
\label{eq:pushforward-mass-measure}
\rmop{pr}_{*}(\mu_{p,M}) = \nu_{p,M}^{\mass}.
\end{equation}
\end{proposition}

\begin{proof}
Recall that we have fixed coordinates, so that \(\bfH(\bbA)\subset \SL_{3}(\bbA)\) and \(\iota(h)=\bbP(h)\).
Since \(\omega_{p}=(\bbP(\xi ),1_{\rmf})\in \bfG(\bbA)\), we have
\[
\iota(h)\omega_{p}=\bbP(h)(\bbP(\xi ),1_{\rmf})=(\bbP(h_{\infty}\xi ),\bbP(h_{\rmf}))=\bbP(a),
\]
where we write \(a=(a_{\infty},a_{\rmf})=(h_{\infty}\xi ,h_{\rmf})\in \GL_{3}(\bbA)\).

We now apply \eqref{equation:projection-shape-lattice-correspondence} for this \(a\).
We have \(N(a_{\rmf})=\Q^{3}\cap h_{\rmf}\hat{\Z}{}^{3}=h_{\rmf}\cdot M\) by \eqref{eq:local-global-lattice}, and \(q_{a_{\infty}}(v)=q_{0}(a_{\infty}^{-1}v)=q_{0}(\xi ^{-1}h_{\infty}^{-1}v)\) with \(q_{0}(w)=w^{T}w\).
Since \(\xi ^{-T}\xi ^{-1}=\Xi^{\circ}\), and \(h_{\infty}\in \bfH(\R)\) satisfies \(h_{\infty}^{T}\Xi^{\circ}h_{\infty}=\Xi^{\circ}\), we get
\begin{align*}
q_{a_{\infty}}(v)
& = (h_{\infty}^{-1}v)^{T}\xi ^{-T}\xi ^{-1}(h_{\infty}^{-1}v)
 = (h_{\infty}^{-1}v)^{T}\Xi^{\circ}(h_{\infty}^{-1}v)\\
& = v^{T}(h_{\infty}^{-T}\Xi^{\circ}h_{\infty}^{-1})v
 = v^{T}\Xi^{\circ}v
 = (\det \Xi)^{-1/3}q(v).
\end{align*}
Then, \eqref{equation:projection-shape-lattice-correspondence} gives
\begin{align*}
\pr(\bfG(\Q)\iota(h)\omega_{p})
&=\pr(\bfG(\Q)\bbP(a))
=\sh(N(a_{\rmf}),q_{a_{\infty}})\\
&=\sh\bigl(h_{\rmf}\cdot M,(\det \Xi)^{-1/3}q\bigr)
=\sh(h_{\rmf}\cdot M,q),
\end{align*}
since the shape of a quadratic lattice is unchanged by scaling.
This proves \eqref{eq:projection-to-shapes}.

The map \(\scrC_{\bfS} \to \scrX\), \([s_{\rmf}]\mapsto \sh(\rho(s_{\rmf})\cdot M,q)\) is well defined, since the isometry class depends only on \([s_{\rmf}]\).
Taking \(h=\rho(s)\) in \eqref{eq:projection-to-shapes} gives, by \eqref{eq:compatibility-homogeneous},
\[
\pr(\Phi_{p}(\bfS(\Q)s))
=\pr(\bfG(\Q)j(s)\omega_{p})
=\pr(\bfG(\Q)\iota\circ \rho(s)\omega_{p})
=\sh(\rho(s_{\rmf})\cdot M,q).
\]
Thus, the diagram
\[
\begin{tikzcd}[column sep=9em]
{}[\bfS] \dar[swap]{\cl} \rar{\Phi_{p}} & X_{\bbA} \dar{\pr}\\
\scrC_{\bfS} \rar[swap]{[s_{\rmf}]\mapsto \sh(\rho(s_{\rmf})\cdot M,q)} & \scrX
\end{tikzcd}
\]
commutes.
By Proposition \ref{proposition:stabilizer-conjugate}, \(\sh(M_{E},q_{E})=\sh(\rho(s_{\rmf})\cdot M,q)\) for \([E]=\Theta^{-1}([s_{\rmf}])\), so the bottom arrow is the composition of \(\Theta^{-1}\) with \([E]\mapsto \sh(M_{E},q_{E})\).
This gives the commutativity of the diagram \eqref{eq:goal-diagram}.

Finally, \(\cl^{-1}(\frakc)=\scrT_{\frakc}\), by Lemma \ref{lemma:fiber-measure}, so \eqref{eq:orbit-volume} gives
\begin{align*}
\cl_{*}m_{\bfS}
& = \sum_{[s_{\rmf}]\in \scrC_{\bfS}} \cl_{*}m_{\bfS}([s_{\rmf}]) \delta_{[s_{\rmf}]}
 = \sum_{[s_{\rmf}]\in \scrC_{\bfS}} m_{\bfS}\bigl(\cl^{-1}([s_{\rmf}])\bigr) \delta_{[s_{\rmf}]}\\
& = \sum_{[s_{\rmf}]\in \scrC_{\bfS}} m_{\bfS}(\scrT_{[s_{\rmf}]}) \delta_{[s_{\rmf}]}
 = \frac{1}{M_{p}}\sum_{[s_{\rmf}]\in \scrC_{\bfS}} \frac{\delta_{[s_{\rmf}]}}{|\calO_{s_{\rmf}}^{\times}|}.
\end{align*}
By commutativity of the diagram, \(\pr\circ\Phi_{p}\) is the composite of \(\cl\) followed by the bottom arrow, so the pushforward \((\pr\circ \Phi_{p})_{*}m_{\bfS}\) is the pushforward of \(\cl_{*}m_{\bfS}\) under the bottom arrow, which sends each \(\delta_{[s_{\rmf}]}\) to \(\delta_{\sh(\rho(s_{\rmf})\cdot M,q)}\).
Hence, by the definitions \eqref{eq:definition-special-homogeneous-set} and \eqref{eq:mass-weighted-measure},
\[
\pr_{*} \mu_{p,M} = (\pr\circ \Phi_{p})_{*}m_{\bfS} = \frac{1}{M_{p}} \sum_{[s_{\rmf}]\in \scrC_{\bfS}} \frac{\delta_{\sh(\rho(s_{\rmf})\cdot M,q)}}{|\calO_{s_{\rmf}}^{\times}|} = \nu_{p,M}^{\mass}.
\]
This completes the proof.
\end{proof}

\section{Effective equidistribution}

We are finally in a position to apply the equidistribution theorem to the homogeneous sets and measures constructed in the previous section.
The remaining analytic step is the comparison of adelic and archimedean Sobolev norms.

\subsection{The adelic equidistribution theorem}

We follow \cite[\S 1.4]{EinsiedlerMargulisMohammadiVenkatesh2020}.
Let \(\bfG(\bbA)^{+}\) denote the image of the simply connected cover, \(\bfG(\bbA)^{+}=\rmop{im}(\bfSL_{3}(\bbA) \to \bfG(\bbA))\), a closed normal subgroup of \(\bfG(\bbA)\).

Recall that \(X_{\bbA}=\bfG(\Q)\bs \bfG(\bbA)\).
Using this we define the decomposition
\[
L^{2}(X_{\bbA}, m_{\bbA})=L^{2}_{0}(X_{\bbA}, m_{\bbA})\oplus L^{2}(X_{\bbA}, m_{\bbA})^{\bfG(\bbA)^{+}},
\]
where \(L^{2}(X_{\bbA}, m_{\bbA})^{\bfG(\bbA)^{+}}\) denotes the space of \(\bfG(\bbA)^{+}\)-invariant functions and \(L^{2}_{0}(X_{\bbA}, m_{\bbA})\) is the orthogonal complement of \(L^{2}(X_{\bbA}, m_{\bbA})^{\bfG(\bbA)^{+}}\).
Note that if \(\bfG\) is simply connected, then \(\bfG(\bbA)=\bfG(\bbA)^{+}\) and \(L^{2}(X_{\bbA}, m_{\bbA})^{\bfG(\bbA)^{+}}\) is the space of constant functions.

The group \(\bfG(\bbA)^{+}\) is a closed, normal subgroup of \(\bfG(\bbA)\); see, e.g., \parencite[p.~451]{PlatonovRapinchuk1994}.
Therefore, the \(L^{2}\)-subspaces in the above decomposition are \(\bfG(\bbA)\)-invariant.
Let
\[
\pi^{+}:L^{2}(X_{\bbA}, m_{\bbA})\to L^{2}(X_{\bbA}, m_{\bbA})^{\bfG(\bbA)^{+}}
\]
denote the orthogonal projection onto the space of \(\bfG(\bbA)^{+}\)-invariant functions.
Let \(C^{\infty}_{c}(X_{\bbA})\) denote the space of smooth compactly supported functions on \(X_{\bbA}\).
Given \(F\in C^{\infty}_{c}(X_{\bbA})\), we define
\[
\pi^{+}F(x)=\int_{X_{\bbA}}F\,d\mu_{x\bfG(\bbA)^{+}}.
\]
So \(\pi^{+}F(x)\) is the average of \(F\) over the closed orbit \(x\bfG(\bbA)^{+}\) with respect to its invariant probability measure.
The main dynamical input is \cite[Theorem 1.5]{EinsiedlerMargulisMohammadiVenkatesh2020}.

\begin{theorem}[EinsiedlerMargulisMohammadiVenkatesh2020]
\label{theorem:EinsiedlerMargulisMohammadiVenkatesh2020}
Let \(\scrD=(\bfT,j_{\bfT},a)\) be a maximal datum as in Definition \ref{definition:period-datum}, with \(\bfT\) simply connected, \(\mu_{\scrD}\) a homogeneous measure.
Then, there exist constants \(C_{0},\kappa_{0}>0\), and a positive integer \(d\), depending only on \((\bfG,\tau, K_{\rmf}, \Omega_{0})\), such that 
\begin{equation}
\label{eq:EinsiedlerMargulisMohammadiVenkatesh2020-bound}
\left|\int_{X_{\bbA}} F\, d\mu_{\scrD} -\pi^{+}F(y)\right| \leq C_{0} (\rmop{vol}(Y_{\scrD}))^{-\kappa_{0}} \calS_{d}(F),
\end{equation}
for every smooth compactly supported \(F\) on \(X_{\bbA}\) and every \(y\in Y_{\scrD}\), where \(\calS_{d}(\cdot)\) is the adelic Sobolev norm \eqref{eq:adelic-Sobolev-norm}.
\end{theorem}

The \(\bfG(\bbA)^{+}\)-invariant term has the following form.

\begin{lemma}
\label{lemma:invariant-term}
For \(f\in C_{c}^{\infty}(\scrX)\) and \(F=f\circ\pr\),
\[
\pi^{+}F(y )= \int_{X_{\bbA}} F \, dm_{\bbA} = \int_{\scrX} f \, dm_{\scrX} \quad(y\in X_{\bbA}).
\]
\end{lemma}

\begin{proof}
By Proposition \ref{proposition:adelic-shape-space}, \(F=f\circ \pr\) is right \(K_{\rmf}\)-invariant.
Then, \cite[(3.2)]{EinsiedlerMargulisMohammadiVenkatesh2020} shows that \(\pi^{+}F=\int_{X_{\bbA}} F \, d m_{\bbA}\).
By \eqref{eq:pushforward-adelic-measure}, \(\pr_{*}m_{\bbA}=m_{\scrX}\) of Proposition \ref{proposition:adelic-shape-space}, so the second equality follows.
\end{proof}

\subsection{Adelic and archimedean Sobolev norms}
\label{subsection:adelic-Sobolev-norm}

In order to apply this result, we need to relate adelic and archimedean Sobolev norms.
In \cite[\S 2.1 and Appendix A]{EinsiedlerMargulisMohammadiVenkatesh2020}, adelic Sobolev norms are defined by using a faithful rational representation, which then defines a collection of compact open subgroups.
Thus, we are not a priori free to choose the compact open subgroup, as it is defined by our correspondence.
Nonetheless, the adelic norm depends only on a limited number of choices, such that it remains independent of \(p\), of \(E_{*}\), and the lattice \(M\).
As a result, for \(F=f\circ \pr\) with \(f\in C_{c}^{\infty}(\scrX)\) as in our case, the adelic norm \(\calS_{d}(F)\) essentially reduces to an archimedean norm of \(f\) on \(\scrY\) comparable to the pre-defined archimedean norm \(\calN_{d}(f)\).
The adelic Sobolev norm depends on the following data:
\begin{enumerate}[auto]
\item The auxiliary representation.
\item The level subgroups.
\item The Lie algebra data.
\item The set of places and the height.
\item The formula for the norm.
\end{enumerate}
We will show that these choices for the adelic Sobolev norm are compatible with our setting for the archimedean Sobolev norm.

\subsection*{(AS1) A faithful representation}
\label{AS1}

Let \(\tau:\bfG \to \bfSL(\frakop{sl}_{3})\cong \bfSL_{8}\) be the conjugation representation, a faithful representation defined by \eqref{eq:faithful-representation}, where \(\frakop{sl}_{3,\Z}\) has the standard Chevalley basis.

\subsection*{(AS2) Finite levels}
\label{AS2}

(\cite[\S 2.1, (2.1)]{EinsiedlerMargulisMohammadiVenkatesh2020}.)
We need to ensure that our compact subgroup \(K_{\rmf}\) is the one defined in \cite[\S 2.1]{EinsiedlerMargulisMohammadiVenkatesh2020}.
The maximal compact subgroup defined in \cite{EinsiedlerMargulisMohammadiVenkatesh2020} is defined from the following finite levels.
For each prime \(\ell\), let
\begin{align*}
K_{\ell}[0]
& = \tau^{-1}\bigl(\SL(\frakop{sl}_{3}(\Z_{\ell}))\bigr),\\
K_{\ell}[m]
& = \rmop{ker}\bigl(K_{\ell}[0] \to \SL(\frakop{sl}_{3}(\Z/\ell^{m}\Z))\bigr)\quad (m\geq 1).
\end{align*}
The maximal compact subgroup of  \cite{EinsiedlerMargulisMohammadiVenkatesh2020} is defined to be \(\prod_{\ell} K_{\ell}[0]\).
We show below that \(K_{\ell}[0]=\bfPGL_{3}(\Z_{\ell})\), for every prime \(\ell\), so that \(\prod_{\ell} K_{\ell}[0]=K_{\rmf}\), that is, our compact open subgroup \(K_{\rmf}\) is the same as the one defined there.

\subsection*{(AS3) Lie algebras and differential operators}
\label{AS3}

(\cite[\S\S 7.1, 7.4, and A.3]{EinsiedlerMargulisMohammadiVenkatesh2020}.)
Let \(\frakg =\rmop{Lie}(\bfG)=\frakop{pgl}_{3,\Q}\) and \(\frakg_{\infty}=\frakg\otimes_{\Q}\R\).
The differential of the quotient map \(\bbP:\bfGL_{3} \to \bfG\) has kernel \(\rmop{Lie}(\bbG_{\rmm})=\Q I_{3}\) and is surjective, so \(\frakg \cong\frakop{gl}_{3}/\Q I_{3}\), and its restriction to \(\frakop{sl}_{3}\), \(d\bbP: \frakop{sl}_{3} \to \frakg \), \(Z\mapsto Z+\Q I_{3}\), is an isomorphism of Lie algebras over \(\Q\).
Using \(d\bbP\), we transport the integral and Euclidean structure of \S \ref{section:setting} to \(\frakg\).
The required bracket stable \(\Z\)-form of \(\frakg\) is \(\frakg_{\Z} = d\bbP(\frakop{sl}_{3}(\Z))\);
stability under the Lie bracket follows from functoriality and since \(\frsl_{3}(\Z)\) is bracket stable.
Note that \(\frakg_{\Z}\) is \emph{not} the Lie algebra of the \(\Z\)-group scheme \(\bfPGL_{3}\), but indeed \(\frakg_{\Z}\otimes_{\Z}\Q\cong \frakg=\frakop{pgl}_{3,\Q}\).
The norm on \(\frakg_{\infty}\) is then \(\|d\bbP(Z)\|_{\frakg,\infty}=\|Z\|_{\frsl_{3}}\) (with the norm on \(\frsl_{3}\) fixed in \S \ref{section:setting}), so that \(d\bbP:\frakop{sl}_{3}(\R)\to \frakg_{\infty}\) is an isometry carrying the basis \(\calB_{\infty}\) to the orthonormal basis \(d\bbP(\calB_{\infty})\).

In \hyperref[AS4]{(AS4)} below, we shall choose \(S=\{\infty\}\), so the norm at the infinite place is precisely the norm \(\|\;\|_{S}\) from \cite[\S 7.1]{EinsiedlerMargulisMohammadiVenkatesh2020}.

Let \(C^{\infty}(X_{\bbA})\) denote the space of functions which are invariant by a compact open subgroup of \(\bfG(\bbA_{\rmf})\).
We shall only consider here a subcollection of these:
let us call a function \(F\in C^{\infty}(X_{\bbA})\) \emph{smooth} if in addition, the function \(\varphi_{x}:h\mapsto F(x(h,1_{\rmf}))\) is smooth on the Lie group \(\bfG(\R)\) for every \(x\in X_{\bbA}\).
Observe that if \(F=f\circ r\) is right \(K_{\rmf}\)-invariant, we have for \(g\in G_{\infty}\), \(\varphi_{x}(\bbP(g))=f(r(x)g)\) by Proposition \ref{proposition:adelic-shape-space}(2).
Then, \(\varphi_{x}\in C^{\infty}(\bfG(\R))\) if and only if \(\varphi_{x}\circ \bbP\in C^{\infty}(\SL_{3}(\R))\) (see \cite[\S\S 2.10, 2.11]{Varadarajan1974}).
so a right \(K_{\rmf}\)-invariant function \(F\) is smooth in our sense if and only if \(F=f\circ r\) with \(f\in C^{\infty}(\scrY)\).
Let \(R_{\bbA}\) be the right regular representation of \(\bfG(\R)\) acting on functions on \(X_{\bbA}\).
For such a function \(F\) and \(Z\in \frsl_{3}(\R)\), and \(x\in X_{\bbA}\), we define
\[
\bigl(dR_{\bbA}(d\bbP(Z))F\bigr)(x) = \left.\frac{d}{dt}\right|_{t=0} F\bigl(x(\exp(td\bbP(Z)),1_{\rmf})\bigr),
\]
where \(\exp:\frakg_{\infty} \to \bfG(\R)\) is the exponential map of the Lie group \(\bfG(\R)\) whose Lie algebra is \(\frakg_{\infty}=d\bbP(\frakop{sl}_{3}(\R))\).
Then, \(dR_{\bbA}(d\bbP(Z))F\) is also smooth;
in fact, \(dR_{\bbA}\) is the derived representation of the right regular representation of \(\bfG(\bbR)\) on smooth functions on \(X_{\bbA}\) through the maps \(\varphi_{x}\).
Then, \(dR_{\bbA}\) is a representation of the Lie algebra \(\frakg_{\infty}\) on \(C^{\infty}(X_{\bbA})\) (see \cite[\S 2.2, (2.2.12)]{Varadarajan1974}), and it extends uniquely to a representation of the universal enveloping algebra \(\frakU(\frakg_{\infty})\) on that space (see \cite[Theorem 3.2.3]{Varadarajan1974}).
Let \(\frakU(d\bbP):\frakU(\frsl_{3}(\R))\to \frakU(\frakg_{\infty})\) be the isomorphism of enveloping algebras induced by \(d\bbP\).
Then for \(D\in \frakU(\frsl_{3}(\R))\), we put \(D_{\bbA}=dR_{\bbA}(\frakU(d\bbP)D)\).
For \(Z\in \frsl_{3}(\R)\), \(Z_{\bbA}=dR_{\bbA}(d\bbP(Z))\).
From \(\bbP:G_{\infty} \to \bfG(\R)\) and naturality of the exponential map, we have \(\exp(t\,d\bbP(Z))=\bbP(\exp(tZ))\) for \(Z\in \frsl_{3}(\R)\), so that
\[
(Z_{\bbA}F)(x) = \left.\frac{d}{dt}\right|_{t=0} F\bigl(x(\bbP(\exp(tZ)),1_{\rmf})\bigr)\quad
(Z\in \frsl_{3}(\R)).
\]
For smooth \(\phi\) on \(\scrY\), we write \(D\phi = dR_{\infty}(D)\phi\) as in \S \ref{section:setting}.

\subsection*{(AS4) The finite set of places and the height.}
\label{AS4}

(See \cite[\S 7.1]{EinsiedlerMargulisMohammadiVenkatesh2020}.)
In the notation of \cite{EinsiedlerMargulisMohammadiVenkatesh2020}, the finite set of places may be taken to be \(S=S_{0}=\{\infty\}\), then \(\Z_{S}=\Z\), \(K^{\prime}(S)=K_{\rmf}\), and the arithmetic subgroup denoted \(\Gamma_{S}=\bfG(\Q)\cap K_{\rmf}\) there is our \(\Delta_{\infty}\), since the base number field is \(\Q\) (Lemma \ref{lemma:class-number-one}).
The space \(X_{S}\) is \(\Delta_{\infty} \bs \bfG(\R)\), and their projection \(\pi_{S}:X_{\bbA} \to X_{S}\) is the quotient map \(X_{\bbA} \to X_{\bbA}/K_{\rmf}\) followed by the (first) homeomorphism \(X_{\bbA}/K_{\rmf}\stosim \Delta_{\infty}\bs \bfG(\R)\) in \eqref{eq:adelic-real-isomorphism}.
The map \(r:X_{\bbA} \to \scrY\) of Proposition \ref{proposition:adelic-shape-space} is the composite of these two maps (i.e., \(\pi_{S}\) of \cite{EinsiedlerMargulisMohammadiVenkatesh2020}) just described, followed by the second homeomorphism \(\Delta_{\infty}\bs \bfG(\R)\to \scrY\), \(\Delta_{\infty}\bbP(g)\mapsto \Gamma_{\infty} g^{\circ}\) of \eqref{eq:adelic-real-isomorphism}.
Every \(h\in \bfG(\R)\) is of the form \(\bbP(g)\) for some \(g\in G_{\infty}\), and for such \(g\in G_{\infty}\), one has \(r(x)=\Gamma_{\infty}g\) (\(x\in X_{\bbA}\)) if and only if \(\pi_{S}(x)=\Delta_{\infty}\bbP(g)\).
Hence, the adelic height defined in \cite[\S 7.1]{EinsiedlerMargulisMohammadiVenkatesh2020}, computed from projective representatives \(h=\bbP(g)\) is as follows.
For \(x\in X_{\bbA}\),
\begin{equation}
\label{eq:adelic-height}
\rmht_{\bbA}(x) 
=\sup\Bigl\{
\|\rmop{Ad}_{\bfG}\bigl(\bbP(g)^{-1}\bigr)(d\bbP(Z))\|_{\frakg,\infty}^{-1} \mid Z\in \frsl_{3}(\Z)\setminus\{0\}\Bigr\},
\end{equation}
where \(r(x)=\Gamma_{\infty}g\), \(g\in G_{\infty}\).
This holds because \(d\bbP(Z)\) runs over \(\frakg_{\Z}\setminus\{0\}\) as \(Z\) runs over \(\frsl_{3}(\Z)\setminus\{0\}\).
This height is independent of the chosen representative \(g\in G_{\infty}\).

\subsection*{(AS5) The formula for the adelic Sobolev norm}
\label{AS5}

(See \cite[Appendix A, \S\S A.3, A.4, (A.3)]{EinsiedlerMargulisMohammadiVenkatesh2020}.)
For a prime \(\ell\) and for \(m\geq 0\), let \(\Av_{\ell}[m]\) be the averaging projection on \(K_{\ell}[m]\)-invariant vectors, i.e., the averaging over \(K_{\ell}[m]\) with respect to Haar probability measure, and put \(\Pi_{\ell}[0]=\Av_{\ell}[0]\) and \(\Pi_{\ell}[m]=\Av_{\ell}[m]-\Av_{\ell}[m-1]\), following \cite[\S A.1]{EinsiedlerMargulisMohammadiVenkatesh2020}, where \(\Pi_{\ell}[\cdot]\) is denoted \(\pr_{\ell}[\cdot]\) there.
Following \cite[\S A.2]{EinsiedlerMargulisMohammadiVenkatesh2020}, we denote by \(\bfm=(m_{\ell})\) a family of nonnegative integers, with \(m_{\ell}=0\) for almost all \(\ell\).
We put \(\|\bfm\|=\prod_{\ell} \ell^{m_{\ell}}\), and put \(\Pi[\bfm]=\prod_{\ell} \Pi_{\ell}[m_{\ell}]\).
Given an integer \(d\geq 0\), the \emph{adelic Sobolev norm} of degree \(d\) is
\begin{equation}
\label{eq:adelic-Sobolev-norm}
\calS_{d}(F)^{2}= \sum_{\bfm} \left(
\|\bfm\|^{d} \sum_{D\in \calD_{d}} \bigl\lVert \Pi[\bfm](1+\rmht_{\bbA})^{d} D_{\bbA}F\bigr\rVert_{L^{2}(X_{\bbA},m_{\bbA})}^{2}
\right).
\end{equation}
(This is \cite[(A.3)]{EinsiedlerMargulisMohammadiVenkatesh2020}.)
Recall that \(\calD_{d}\subset \frakU(\frsl_{3}(\R))\) was introduced in \S \ref{section:setting}, and there is only one archimedean place.

We now verify the compatibility of our choices with the adelic and archimedean Sobolev norms.

\begin{lemma}
\label{lemma:adelic-archimedean-equidistribution-objects}
\begin{enumerate}[auto]
\item For every \(a\in \bfG(\R)\) and every \(Z\in \frakop{sl}_{3}(\R)\),
\[
\rmop{Ad}_{\bfG}(a) d\bbP(Z)=d\bbP\bigl(\tau(a)Z\bigr).
\]
\item We have \(K_{\ell}[0]=\bfPGL_{3}(\Z_{\ell})\) for every prime \(\ell\), hence \(\prod_{\ell} K_{\ell}[0]=K_{\rmf}\), and every \(K_{\ell}[m]\subset K_{\rmf}\).
\item The adelic and archimedean heights satisfy \(\rmht_{\bbA} = \rmht_{\infty}\circ r\).
In particular, \(\rmht_{\bbA}\) is right \(K_{\rmf}\)-invariant.
\item For every \(\phi\in C^{\infty}(\scrY)\) and every \(D\in \frakU(\frakop{sl}_{3}(\R))\), we have \(D_{\bbA}(\phi\circ r)=(D\phi)\circ r\).
In particular, for \(f\in C_{c}^{\infty}(\scrX)\), we have \(F=f\circ \pr\) and \(\hat{f}(y)=f(yK_{\infty})\), and
\[
D_{\bbA} F=(D\hat{f})\circ r.
\]
\end{enumerate}
\end{lemma}

\begin{proof}
(1) From \(a=\bbP(g)\) for some \(g\in \GL_{3}(\R)\), we have \(\tau(a)Z=gZg^{-1}\), and since the differential \(d\bbP\) is \(\rmop{Ad}\)-equivariant (see \cite[10.21, (61)]{Milne2017}), we get \(d\bbP(gZg^{-1})=\rmop{Ad}_{\bfG}(\bbP(g))d\bbP(Z)\).
(In fact, this is functorial, so that \(\tau=(d\bbP)^{-1}\circ \Ad_{\bfG}\circ d\bbP\) as \(\Q\)-morphisms.)

(2) Conjugation by an element \(g\in \GL_{3}(\Z_{\ell})\) preserves \(\frsl_{3}(\Z_{\ell})\), i.e., \(\tau(\bbP(g))\in \SL(\frsl_{3}(\Z_{\ell}))\).
Since \(K_{\ell}[0]=\tau^{-1}(\SL(\frsl_{3}(\Z_{\ell})))\), this means \(\bbP(g)\in K_{\ell}[0]\), hence \(\bfPGL_{3}(\Z_{\ell})\subset K_{\ell}[0]\).
For the reverse inclusion, for \(a\in K_{\ell}[0]\), we use the Cartan decomposition to write a lift of \(a\) in \(\GL_{3}(\Q_{\ell})\) as \(u\delta v\) with \(\delta=\rmop{diag}(\ell^{a_{1}},\ell^{a_{2}},\ell^{a_{3}})\) and \(u,v\in \GL_{3}(\Z_{\ell})\), \(a_{i}\in \Z\).
Since conjugation by \(u\) and \(v\) preserve \(\frsl_{3}(\Z_{\ell})\), conjugation by \(\delta\) must too.
Conjugation by \(\delta\) maps the matrix unit \(E_{ij}\) to \(\ell^{a_{i}-a_{j}}E_{ij}\), hence we must have \(a_{i}-a_{j}\geq 0\) for all \(i\neq j\), i.e., \(a_{1}=a_{2}=a_{3}\), and \(u\delta v \in \Q_{\ell}^{\times} \GL_{3}(\Z_{\ell})\), hence \(a\in \bfPGL_{3}(\Z_{\ell})\).
Thus, \(\prod_{\ell} K_{\ell}[0]=\prod_{\ell} \bfPGL_{3}(\Z_{\ell})=K_{\rmf}\) by definition.

(3) Let \(x\in X_{\bbA}\) and \(r(x)=\Gamma_{\infty}g\) with \(g=g^{\circ}\in G_{\infty}\) as in Proposition \ref{proposition:adelic-shape-space}.
For \(Z\in \frsl_{3}(\Z)\), part (1) with \(a=\bbP(g)^{-1}=\bbP(g^{-1})\) gives
\[
\Ad_{\bfG}(\bbP(g)^{-1})d\bbP(Z)=d\bbP(g^{-1}Zg)=d\bbP\bigl(\Ad(g^{-1})Z\bigr),
\]
Thus the norm \(\|\,\|_{\frakg,\infty}\) in the definition of the adelic height \eqref{eq:adelic-height} is \(\|d\bbP(\Ad(g^{-1})Z)\|_{\frakg,\infty}=\|\Ad(g^{-1})Z\|_{\frsl_{3}}\) since the norm on \(\frakg_{\infty}\) is obtained by transport from the norm on \(\frsl_{3}(\R)\).
Hence, the supremum in the adelic height \eqref{eq:adelic-height} coincides with the supremum in the archimedean height \eqref{eq:archimedean-height} at \(\Gamma_{\infty}g\), i.e., \(\rmht_{\bbA}(x)=\rmht_{\infty}(\Gamma_{\infty}g)=\rmht_{\infty}(r(x))\).
Right \(K_{\rmf}\)-invariance is immediate from right \(K_{\rmf}\)-invariance of \(r\) by Proposition \ref{proposition:adelic-shape-space}.

(4) Let \(Z\in \frsl_{3}(\R)\) and \(x\in X_{\bbA}\).
Then, \(\exp(tZ)\in G_{\infty}\) so \(\exp(tZ)=\exp(tZ)^{\circ}\).
Therefore, Proposition \ref{proposition:adelic-shape-space}(2) gives
\[
r(x(\bbP(\exp(tZ)),1_{\rmf}))=r(x)(\exp(tZ))^{\circ}=r(x)\exp(tZ).
\]
Together with the formula for \(Z_{\bbA}\) in \hyperref[AS3]{(AS3)}, we get for \(\phi\in C^{\infty}(\scrY)\),
\[
Z_{\bbA}(\phi \circ r)(x)=\left.\frac{d}{dt}\right|_{t=0} \phi\bigl(r(x)\exp(tZ)\bigr) 
= \bigl(dR_{\infty}(Z)\phi\bigr)(r(x)),
\]
by the definition of \(dR_{\infty}(Z)\) in \S \ref{section:setting}.
Thus,
\[
Z_{\bbA}(\phi \circ r)=(dR_{\infty}(Z)\phi)\circ r=(Z\phi)\circ r.
\]
This is the desired assertion for the special case \(D=Z\) for every smooth \(\phi\).
The collection of operators \(D\in \frakU(\frakop{sl}_{3}(\R))\) such that \(D_{\bbA}(\phi \circ r) = (D\phi)\circ r\) for every \(\phi\in C^{\infty}(\scrY)\) is a linear subspace, clearly closed under products, hence a subalgebra;
it contains \(\frakop{sl}_{3}(\R)\), which generates the enveloping algebra.
It follows that the identity \(D_{\bbA}(\phi\circ r)=(D\phi)\circ r\) holds for all \(D\in \frakU(\frsl_{3}(\R))\).
The last claim follows with \(\phi=\hat{f}\) since \(F=f\circ \pr = \hat{f} \circ r\) from Proposition \ref{proposition:adelic-shape-space}(6).
\end{proof}

\subsection{Comparison of the Sobolev norms}

We now show that we can compare the archimedean and adelic Sobolev norms.
For this, we will need a lower bound on the archimedean height.

\begin{lemma}
One has \(c_{h} = \inf_{y\in \scrY} \rmht_{\infty}(y)>0\).
\end{lemma}

Note that this is independent of \(p\), \(E_{*}\), or \(M\).
We put \(C_{h}=1+c_{h}^{-1}\).

\begin{proof}
Let \(g\in G_{\infty}\).
The automorphism \(\Ad(g^{-1})\) of \(\frsl_{3}(\R)\) has determinant one, so \(\Ad(g^{-1})\frsl_{3}(\Z)\) is a lattice in the Euclidean space \(\frsl_{3}(\R)\) equipped with the norm \(\|\;\|\), with the same covolume as \(\frsl_{3}(\Z)\).
Let \(c_{0}\) be this volume.
By Minkowski's convex body theorem (\cite{Cassels1959}), there exists a radius \(R>0\) depending only on \(c_{0}\) and the norm \(\|\;\|\), such that every such lattice \(\Ad(g^{-1})\frsl_{3}(\Z)\), since it has covolume \(c_{0}\), contains a nonzero vector of norm at most \(R\), i.e., for every \(g\in G_{\infty}\), there is \(Z_{g}\in \frsl_{3}(\Z)\setminus\{0\}\) such that \(\|\Ad(g^{-1})Z_{g}\|\leq R\), i.e., \(\rmht_{\infty}(\Gamma_{\infty}g)\geq 1/R\) for every \(g\in G_{\infty}\), so \(c_{h}\geq 1/R>0\).
\end{proof}

\begin{proposition}
\label{proposition:Sobolev-norms-comparison}
\begin{enumerate}[auto]
\item Let \(d\geq 0\), \(f\in C_{c}^{\infty}(\scrX)\), and \(\hat{f}(y)=f(yK_{\infty})\).
Then,
\begin{equation}
\label{eq:adelic-Sobolev-norm-simple}
\calS_{d}(f\circ \pr)^{2} = \sum_{D\in \calD_{d}} \int_{\scrY} (1+\rmht_{\infty})^{2d} |D\hat{f}|^{2} \,dm_{\scrY},
\end{equation}
and
\begin{equation}
\label{eq:adelic-Sobolev-norm-comparison}
\calN_{d}(f)\leq \calS_{d}(f\circ \pr) \leq C_{h}^{d} \calN_{d}(f).
\end{equation}
\item There exists an integer \(d_{0}\geq 1\), and for every \(d\geq d_{0}\), there exists a constant \(C_{0,d}\) depending only on our data, such that
\[
\|f\|_{\infty} \leq C_{0,d} \, \calN_{d}(f)\quad
(f\in C_{c}^{\infty}(\scrX)).
\]
\end{enumerate}
\end{proposition}

\begin{proof}
(1) For \(D\in \calD_{d}\), let \(Q_{D}=(1+\rmht_{\bbA})^{d} D_{\bbA}F\).
By Lemma \ref{lemma:adelic-archimedean-equidistribution-objects}, (3)--(4), \(\rmht_{\bbA}=\rmht_{\infty}\circ r\) and \(D_{\bbA}F=(D\hat{f})\circ r\), so
\[
Q_{D} = \psi_{D}\circ r,\quad
\psi_{D}=(1+\rmht_{\infty})^{d}D\hat{f}~\text{on}~\scrY,
\]
and \(Q_{D}\) is right \(K_{\rmf}\)-invariant.

Only \(\bfm=0\) contributes to \eqref{eq:adelic-Sobolev-norm}.
Since \(K_{\ell}[m]\subset K_{\ell}[0]\subset K_{\rmf}\), and \(Q_{D}\) is right \(K_{\rmf}\)-invariant, we have \((\Av_{\ell}[m]Q_{D})(x)=\int_{K_{\ell}[m]} Q_{D}(xk)dk=Q_{D}(x)\), so \(\Av_{\ell}[m]Q_{D}=Q_{D}\) for all \(\ell\) and all \(m\geq 0\).
So \(\Pi[0]Q_{D}=Q_{D}\) and for \(\bfm\neq 0\), \(\Pi[\bfm]Q_{D}=0\).
Thus, \eqref{eq:adelic-Sobolev-norm} only has one term for \(\|\bfm\|=1\), so it becomes
\[
\calS_{d}(F)^{2} = \sum_{D\in \calD_{d}} \|Q_{D}\|_{L^{2}(X_{\bbA},m_{\bbA})}^{2}.
\]

Since \(r_{*}m_{\bbA}=m_{\scrY}\),
\[
\|Q_{D}\|_{L^{2}(X_{\bbA},m_{\bbA})}^{2} =
\|\psi_{D}\circ r\|_{L^{2}(X_{\bbA},m_{\bbA})}^{2} 
= \int_{\scrY} (1+\rmht_{\infty})^{2d}|D\hat{f}|^{2} \, dm_{\scrY}.
\]
Thus, \eqref{eq:adelic-Sobolev-norm-simple} follows by summing over \(D\in \calD_{d}\).
Since \(\rmht_{\infty}\geq c_{h}\), we have \(\rmht_{\infty}\leq 1+\rmht_{\infty}\leq c_{h}^{-1}\rmht_{\infty}+\rmht_{\infty}=(1+c_{h}^{-1})\rmht_{\infty}=C_{h}\rmht_{\infty}\).
Thus, \(\rmht_{\infty}\leq 1+\rmht_{\infty}\leq C_{h} \rmht_{\infty}\).
Together with the definition of the archimedean Sobolev norm from \eqref{eq:real-Sobolev-norm}, this immediately implies the inequalities of \eqref{eq:adelic-Sobolev-norm-comparison}.

(2) By the Sobolev embedding property (see \cite[\S 7.5, Property S1]{EinsiedlerMargulisMohammadiVenkatesh2020}), there exists \(d_{0}\) depending only on \(\dim\bfG\) such that for all \(d\geq d_{0}\), we have \(\|F\|_{L^{\infty}} \leq C_{d} \calS_{d}(F)\) for every smooth compactly supported function \(F\) on \(X_{\bbA}\).
Applied to \(F=f\circ \pr\), since \(\|f\|_{\infty}=\|F\|_{\infty}\), we get
\[
\|f\|_{\infty} \leq C_{d} \calS_{d}(F) \leq C_{d} C_{h}^{d} \calN_{d}(f).
\]
This completes the proof.
\end{proof}

We control the passage from the mass weights \(|\Aut(E)|^{-1}\) to uniform sampling as follows.

\begin{lemma}
\label{lemma:uniform-vs-mass-weighted-measures}
Every function \(\phi:\rmSS_{p}\to \C\) satisfies
\[
\left|\frac{1}{M_{p}} \sum_{[E]} \frac{\phi(E)}{|\Aut(E)|} - \frac{1}{h_{p}} \sum_{[E]} \phi(E)\right| \leq \frac{28}{p-1}\|\phi\|_{\infty}.
\]
\end{lemma}

\begin{proof}
Let \(I\neq\emptyset\) be a finite set, and let \(0<w_{i} \leq b\) for \(i\in I\) be real numbers.
Put \(w^{*}=\sum_{i} w_{i}\) and \(\delta=\sum_{i} (b-w_{i})\).
Then, every complex-valued function \(\phi:I \to \C\) satisfies
\[
\left|\frac{1}{w^{*}} \sum_{i} w_{i} \phi(i) - \frac{1}{|I|} \sum_{i} \phi(i)\right|\leq \frac{2\delta}{w^{*}} \|\phi\|_{\infty},
\]
as is easily checked by an elementary computation.
We apply this with \(I=\rmSS_{p}\), \(w_{[E]}=|\Aut(E)|^{-1}\) and \(b=1/2\).
So \(|I|=h_{p}\) and \(w^{*}=\sum_{[E]} w_{[E]}=M_{p}=\frac{p-1}{24}\).
For \(p>3\), \(|\Aut(E)|=2\) except for two \(j\)-invariants (see \cite{Silverman2009}).
So we can take \(b=1/2\), and except for these two cases, \(b-|\Aut(E)|^{-1}=0\).
For the cases \(j(E)=0\) and \(j(E)=1728\), we have automorphism groups of size \(6\) and \(4\) respectively, so \(b-|\Aut(E)|^{-1}\) is \(1/3\) and \(1/4\) respectively.
Thus, \(\delta=\sum_{[E]} r_{E}\leq 1/3+1/4=7/12\).
If we write \(r_{E}=b-w_{[E]}=(1/2)-w_{[E]}\), we have
\[
M_{p} = \sum_{[E]} w_{[E]} = \sum_{[E]} \bigl(\tfrac{1}{2}-r_{E}\bigr)
=\frac{h_{p}}{2}-\delta.
\]
Since \(|I|=h_{p}\), and \(w^{*}=M_{p}\), the inequality then reads
\[
\left|\frac{1}{M_{p}} \sum_{[E]} \frac{\phi(E)}{|\Aut(E)|} - \frac{1}{h_{p}} \sum_{[E]} \phi(E)\right| \leq \frac{2\cdot 7}{12}\frac{24}{(p-1)}\|\phi\|_{\infty}=\frac{28}{p-1}\|\phi\|_{\infty}.\qedhere
\]
\end{proof}

\subsection{The main equidistribution theorem}

We can now prove our main theorem, which we reformulate here.
Let \(M\in\{L,\Lambda\}\), and fix the basepoint \(y_{p,M}=\bfG(\Q)\omega_{p,M}\in Y_{p,M}\).
By \cite[(A.4)]{EinsiedlerMargulisMohammadiVenkatesh2020}, the adelic Sobolev norm satisfies \(\calS_{d}(F)\leq \calS_{d^{\prime}}(F)\) if \(d<d^{\prime}\), so we may consider a larger \(d\) to ensure that \(d\geq d_{0}\).

\begin{theorem}
\label{theorem:main-equidistribution}
There exists an integer \(d\geq d_{0}\), and there exist constants \(C,\kappa>0\), independent of \(p\), such that the following hold for \(M\in \{L,\Lambda\}\).
\begin{enumerate}[auto]
\item For every \(f\in C_{c}^{\infty} (\scrX)\) and every prime \(p\geq 7\),
\begin{equation}
\label{eq:equidistribution-bound}
\left| \int_{\scrX}f\,d\nu_{p,M}^{\mass}- \int_{\scrX}f\,dm_{\scrX}\right| \leq C p^{-\kappa} \calN_{d}(f).
\end{equation}
\item For every \(f\in C_{c}^{\infty} (\scrX)\) and every prime \(p\geq 7\),
\[
\left| \int_{\scrX}f\,d\nu_{p,M}- \int_{\scrX}f\,dm_{\scrX}\right| \leq C p^{-\kappa} \calN_{d}(f).
\]
\item \(\nu_{p,M}^{\mass}\) and \(\nu_{p,M}\) converge weakly to \(m_{\scrX}\) as \(p\to\infty\) through primes.
\end{enumerate}
\end{theorem}

\begin{proof}
(1): Let \(F=f\circ \rmop{pr}\).
\begin{itemize}[leftmargin=1em]
\item By Proposition \ref{proposition:structure-homogeneous-set}, \(\scrD_{p,M}\) is a maximal datum, hence satisfies the assumptions of Theorem \ref{theorem:EinsiedlerMargulisMohammadiVenkatesh2020}.
\item By Proposition \ref{proposition:homogeneous-set-volume}, equation \eqref{eq:volume-homogeneous-set}, the volume of the algebraic homogeneous set is \(v_{p}=\rmop{vol}(Y_{p,M})=(p-1)/12\).
\item By \eqref{eq:pushforward-mass-measure} of Proposition \ref{proposition:final-projection}, \(\pr_{*}\mu_{p,M}=\nu_{p,M}^{\mass}\).
\item By Lemma \ref{lemma:invariant-term}, \(\pi^{+}F(y)=\int_{\scrX} f\, dm_{\scrX}\) for every \(y\in X_{\bbA}\).
\item Finally, the adelic and archimedean Sobolev norms are related by \eqref{eq:adelic-Sobolev-norm-comparison} from Proposition \ref{proposition:Sobolev-norms-comparison}.
\end{itemize}
Thus, an application of Theorem \ref{theorem:EinsiedlerMargulisMohammadiVenkatesh2020} with \(y=y_{p,M}\) yields
\begin{align*}
\left| \int_{\scrX}f\,d\nu_{p,M}^{\mass}- \int_{\scrX}f\,dm_{\scrX}\right|
&=\left| \int_{X_{\bbA}}F\,d\mu_{p,M}- \pi^{+}F(y_{p,M})\right|\\
&\leq C_{0} \left(\frac{12}{p-1}\right)^{\kappa_{0}} \calS_{d}(F)\\
&\leq C_{0} \left(\frac{12}{p-1}\right)^{\kappa_{0}} C_{h}^{d}\calN_{d}(f).
\end{align*}
For \(p\geq 7\), we have \(p-1\geq 6p/7\), so \eqref{eq:equidistribution-bound} holds with \(\kappa=\kappa_{0}\) and \(C=C_{0}14^{\kappa_{0}}C_{h}^{d}\), which proves (1).

(2)
From the definitions \eqref{eq:uniform} and \eqref{eq:mass} of the measures \(\nu_{p,M}\) and \(\nu_{p,M}^{\mass}\), Lemma \ref{lemma:uniform-vs-mass-weighted-measures} applied with \(\phi(E)=f(\sh(M_{E},q_{E}))\) now reads
\[
\left|\int_{\scrX} f \, d\nu_{p,M}^{\mass} - \int_{\scrX} f \, d\nu_{p,M}\right|
\leq \frac{28}{p-1}\|f\|_{\infty}.
\]
By Proposition \ref{proposition:Sobolev-norms-comparison}(2), we get
\[
\left|\int_{\scrX} f \, d\nu_{p,M}^{\mass} - \int_{\scrX} f \, d\nu_{p,M}\right|
\leq \frac{28}{p-1}C_{0,d}\calN_{d}(f).
\]
Combining this with \eqref{eq:equidistribution-bound} and the triangle inequality, we obtain (2) with \(\kappa=\min\{\kappa_{0},1\}\).

(3) The argument is similar to the proof of \cite[Corollary in \S 1.7]{EinsiedlerMargulisMohammadiVenkatesh2020} (the proof is in \cite[\S 7.12]{EinsiedlerMargulisMohammadiVenkatesh2020}).
Let \(\pi_{K}:\scrY\to \scrX\) be the quotient map.
Let \(C\subset \scrX\) be a compact set, and \(U\supset C\) open.
Choose some \(\chi\in C_{c}^{\infty}(\scrX)\) with \(1_{C}\leq \chi\leq 1_{U}\).
For instance, from the standard construction in \cite[Chapter I, Lemma 1.5]{Helgason2001}, we may pick a function \(\psi\in C_{c}^{\infty}(\scrY)\) such that \(0\leq \psi\leq 1\), and \(1_{\pi_{K}^{-1}(C)}\leq \psi \leq 1_{\pi_{K}^{-1}(U)}\), i.e., \(\psi=1\) on \(\pi_{K}^{-1}(C)\) and supported on \(\pi_{K}^{-1}(U)\).
Then, the average \(y\mapsto \int_{K_{\infty}} \psi(yk)\, dk\) over \(K_{\infty}\) is right \(K_{\infty}\)-invariant and is the lift of the desired map \(\chi\).

Now, let \(\nu\) denote either of the two measures.
Let \(U\subset \scrX\) be an open set, and let \(\epsilon>0\) be arbitrary.
Since \(\scrX\) is locally compact, \(U\) is an increasing union of compact sets, so by inner regularity of the measure \(m_{\scrX}\), there is a compact set \(C\subset U\) such that \(m_{\scrX}(C)> m_{\scrX}(U)-\epsilon\).
Choose \(\chi\in C_{c}^{\infty}(\scrX)\) such that \(1_{C}\leq \chi\leq 1_{U}\).
Integrating, we get \(\nu(C)\leq\nu(\chi)\leq \nu(U)\).
Since \(\nu(\chi)\to m_{\scrX}(\chi)\) as \(p\to\infty\) for both choices of \(\nu\) by parts (1) and (2), we get
\[
\liminf_{p\to\infty} \nu(U) \geq 
\lim_{p\to\infty} \nu(\chi) = m_{\scrX}(\chi) \geq m_{\scrX}(C)>m_{\scrX}(U)-\epsilon.
\]
Letting \(\epsilon\downarrow 0\) shows that \(\liminf_{p\to\infty} \nu(U) \geq m_{\scrX}(U)\) for every open \(U\subset \scrX\).
Since \(\scrX\) is metrizable, this is equivalent to weak convergence of \(\nu\) to \(m_{\scrX}\) by the Portmanteau theorem (e.g., \cite[Theorem 2.1]{Billingsley1999}).
\end{proof}

\section{Successive minima}

For a positive definite quadratic lattice \((M,q)\) of rank \(n\), we let \(\lambda_{i}(M,q)\) denote the successive minima (see \eqref{eq:successive-minima-definition}).
Let
\[
\calW = \{x\in \R_{>0}^{3}\mid x_{1} \leq x_{2} \leq x_{3} \},
\]
and let \(\bflambda:\scrX\to\calW\) be defined by
\[
\bflambda(\Gamma_{\infty} g K_{\infty}) = \bigl(\lambda_{i}(g^{-1}\Z^{3},q_{0})\bigr)_{i=1}^{3}.
\]
Finally, let \(\eta=\bflambda_{*} m_{\scrX}\), so that \(\eta\) is a probability measure on \(\R_{>0}^{3}\) supported on \(\calW\).

\begin{lemma}
The map \(\bflambda:\scrX\to \calW\) is continuous and \(\eta(\{x\mid x_{i}=t\})=0\) for every \(i=1,2,3\) and \(t>0\).
\end{lemma}

\begin{proof}
Continuity of \(\bflambda\) is standard, see e.g., \cite{PettetSouto2008}.
This is easily proved as follows.
For \(u,g\in G_{\infty}\), we have
\[
\|u\|_{\text{op}}^{-1}\, \lambda_{i}(g^{-1}\Z^{3})\leq
\lambda_{i}\bigl((gu)^{-1}\Z^{3}\bigr)\leq
\|u^{-1}\|_{\text{op}}\, \lambda_{i}(g^{-1}\Z^{3}).
\]
Letting \(u\to 1\) shows that \(\bflambda\) is continuous at \(g\), so on \(G_{\infty}\), and since \(\pi_{K}\) is a quotient map, \(\bflambda\) is continuous on \(\scrX\).

For \(y=\Gamma_{\infty}g\), the lattice \(g^{-1}\Z^{3}\) only depends on \(y\), since \((\gamma g)^{-1}\Z^{3}=g^{-1}\Z^{3}\) for \(\gamma\in \Gamma_{\infty}\).
This identifies \(\scrY\) with the space of unimodular lattices in \(\R^{3}\).
The second claim follows immediately from Siegel's mean value theorem (see \cite[Eq. (2), Theorem]{Siegel1945}), which states that for a compactly supported bounded Borel measurable function \(\psi:\R^{3}\to \R\), we have
\[
\int_{\scrY} \sum_{x\in g^{-1}\Z^{3}\setminus\{0\}} \psi(x) \, dm_{\scrY}(\Gamma_{\infty}g) = \int_{\R^{3}} \psi(x)\, dx.
\]
Let \(S_{t}=\{x\in \R^{3}\mid q_{0}(x)=t^{2}\}\) be the sphere of radius \(t\).

Since the minima are attained, \(\lambda_{i}(g^{-1}\Z^{3},q_{0})=t\) implies that \(g^{-1}\Z^{3}\cap S_{t}\neq\emptyset\).
We then apply Siegel's mean value theorem to the indicator function \(1_{S_{t}}\).
Since \(\eta=\bflambda_{*}m_{\scrX}\) and \(m_{\scrX}\) is the pushforward of \(m_{\scrY}\) under the quotient map \(\scrY\to \scrY/K_{\infty}\),
\begin{align*}
\eta(\{x\mid x_{i}=t\})
&= m_{\scrX}\bigl(\bflambda^{-1}(\{x\mid x_{i}=t\})\bigr)\\
&= m_{\scrY}\bigl(\{\Gamma_{\infty}g\in \scrY\mid \lambda_{i}(g^{-1}\Z^{3},q_{0})=t\}\bigr)\\
&\leq \int_{\scrY} |(g^{-1}\Z^{3}\setminus \{0\})\cap S_{t}| \,dm_{\scrY}(\Gamma_{\infty}g) = \rmop{vol}(S_{t})=0.\qedhere
\end{align*}
\end{proof}

We rescale the quadratic form \(q_{E}\) to renormalize the determinants.
For \(M\in \{L,\Lambda\}\), define \(a_{M,p}>0\) by
\begin{equation}
\label{eq:rescaling-factors}
a_{L,p}=(p/2)^{1/3},\quad
a_{\Lambda,p}=(2p)^{1/3}.
\end{equation}
Dividing \(q_{E}\) by \(a_{M,p}^{2}\) normalizes the Gram determinants to \(1\).
Thus,
\[
\bflambda\bigl(\sh(M_{E},q_{E})\bigr)
=\bigl(\lambda_{i}(M_{E},q_{E})/a_{M,p}\bigr)_{1\leq i\leq 3} = u_{p,E}^{(M)}.
\]
Hence, we have \(\eta_{p}^{(M)}=\bflambda_{*} \nu_{p,M}\) and \(\eta_{p}^{(M,\mass)}=\bflambda_{*} \nu_{p,M}^{\mass}\).

We can now prove Theorem \ref{theorem:B}.

\begin{theorem}
\begin{enumerate}[auto]
\item As \(p\to\infty\) through the primes, for either choice \(M\in\{L,\Lambda\}\),
\begin{equation}
\label{eq:weak-convergence-empirical}
\eta_{p}^{(M)}\Longrightarrow \eta,\quad
\eta_{p}^{(M,\mass)} \Longrightarrow \eta.
\end{equation}
\item For every Borel set \(\calR\subset \R_{>0}^{3}\) with \(\eta(\partial\calR)=0\), and for \(M\in\{L,\Lambda\}\), we have
\[
\lim_{p\to\infty} \frac{\#\{[E]\in \rmSS_{p} \mid(\lambda_{i}(M_{E},q_{E})/a_{M,p})_{1\leq i\leq 3}\in \calR\}}{h_{p}} = \eta(\calR).
\]
\item For every function \(\psi(p)>0\) such that \(\psi(p)\to\infty\), the proportion of classes of curves satisfying simultaneously for all \(i=1,2,3\),
\begin{equation}
\label{eq:balanced}
\frac{a_{M,p}}{\psi(p)}\leq \lambda_{i}(M_{E},q_{E})\leq a_{M,p}\psi(p),\quad
\frac{\lambda_{3}(M_{E},q_{E})}{\lambda_{1}(M_{E},q_{E})} \leq \psi(p),
\end{equation}
tends to \(1\).
The same results hold for mass-weighted measures with weights \(1/(M_{p}|\Aut(E)|)\) instead of \(1/h_{p}\).
\end{enumerate}
\end{theorem}

\begin{proof}
(1)
The measures are precisely \(\bflambda_{*}\nu_{p,M}\) and \(\bflambda_{*}\nu_{p,M}^{\mass}\).
Since \(\bflambda\) is continuous, by Part (3) of Theorem \ref{theorem:main-equidistribution} and \cite[Theorem 2.7]{Billingsley1999}, the result follows.

(2) The equivalent statement for Borel sets \(\calR\) with \(\eta(\partial \calR)=0\) is from \cite[Theorem 2.1(v)]{Billingsley1999}.

(3)
Let \(P_{p}\) denote either the uniform probability measure on \(\rmSS_{p}\), with measure \(1/h_{p}\) for each class, or the normalized mass measure, with measure \(1/(M_{p}|\Aut(E)|)\) for the class \([E]\).
Thus, for \(A\subset \rmSS_{p}\), we have
\[
P_{p}(A) = \frac{|A|}{h_{p}},\quad
P_{p}(A) = \frac{1}{M_{p}} \sum_{[E]\in A} \frac{1}{|\Aut(E)|}.
\]
Let \(\eta_{p,M}\) denote the pushforward of \(P_{p}\) under the map \(u_{p}^{(M)}:[E]\mapsto u_{p,E}^{(M)}\), that is, \(\eta_{p,M}\) is either \(\eta_{p}^{(M)}\) or \(\eta_{p}^{(M,\mass)}\).
By \eqref{eq:weak-convergence-empirical}, \(\eta_{p,M}\Rightarrow \eta\).
We have observed that \(u_{p,E}^{(M)}(i)= \lambda_{i}(M_{E},q_{E})/a_{M,p}\).

For \(R>1\), consider the open box \(U_{R}=(1/R,R)^{3}\subset\R_{>0}^{3}\), \(R>1\).
For \(R>1\), consider the sets
\[
A_{p}^{(M)}(R)=\{[E]\in \rmSS_{p} \mid u_{p,E}^{(M)}\in U_{R}\},\quad
A_{p}(R) = A_{p}^{(L)}(R)\cap A_{p}^{(\Lambda)}(R).
\]
As \(R\to\infty\), the sets \(U_{R}\) cover \(\R_{>0}^{3}\), so \(\eta(U_{R})\to 1\) as \(R\to\infty\).
Let \(\epsilon>0\) and choose \(R\) such that \(\eta(U_{R})>1-\epsilon\).
Since \(\eta(\partial U_{R})=0\), by the Portmanteau theorem \cite[Theorem 2.1]{Billingsley1999}, we have
\[
\liminf_{p\to\infty} P_{p}(A_{p}^{(M)}(R))
=\liminf_{p\to\infty} \eta_{p,M}(U_{R})\geq \eta(U_{R}) > 1-\epsilon,
\]
for either \(M\in \{L,\Lambda\}\).
Thus, there exists \(p_{0}\) such that both probabilities are \(>1-\epsilon\) for all \(p\geq p_{0}\).
Then, \(P_{p}(\rmSS_{p}\setminus A_{p}(R))\leq P_{p}(\rmSS_{p}\setminus A_{p}^{(L)}(R))+P_{p}(\rmSS_{p}\setminus A_{p}^{(\Lambda)}(R))< 2\epsilon\).
Thus, \(P_{p}(A_{p}(R))> 1-2\epsilon\).

Let \(B_{p}\subset\rmSS_{p}\) denote the set of classes satisfying the inequalities  \eqref{eq:balanced}.
Since \(\psi(p)\to \infty\), with \(R\) fixed, we may increase \(p_{0}\) to \(p_{1}\) such that \(\psi(p)\geq R^{2}\) for every prime \(p\geq p_{1}\).
Let \(p\geq p_{1}\), \([E]\in A_{p}(R)\), \(M\in\{L,\Lambda\}\) and \(i\in \{1,2,3\}\).
Write \(\lambda_{i}=\lambda_{i}(M_{E},q_{E})\) and \(a=a_{M,p}\) for simplicity, so that \(u=(\lambda_{i}/a)_{i}\).
For \([E]\in A_{p}(R)\), we have \(u\in U_{R}\).
Then, \(R^{-1}<\lambda_{i}/a<R\) for each \(i\), and since \(1<R<R^{2}\leq \psi(p)\), we then have
\[
\frac{a}{\psi(p)}\leq \frac{a}{R}<\lambda_{i} <aR< a \psi(p),
\]
and
\[
\frac{\lambda_{3}}{\lambda_{1}} = \frac{\lambda_{3}/a}{\lambda_{1}/a} < \frac{R}{R^{-1}}=R^{2} \leq \psi(p).
\]
Thus, the inequalities \eqref{eq:balanced} hold, and we conclude that \(A_{p}(R)\subset B_{p}\).
Thus, \(P_{p}(B_{p})>1-2\epsilon\) for all primes \(p\geq p_{1}\).
Since \(\epsilon>0\) was arbitrary, \(P_{p}(B_{p}) \to 1\), which proves the assertions for both the uniform measure and the mass-weighted measure.
\end{proof}

\printbibliography

\end{document}